\documentclass[12pt]{article}
\usepackage{amsfonts}
\usepackage{amssymb}

\input{tcilatex}
\begin{document}

\begin{center}
{\Large Beyond Haar and Cameron-Martin: the Steinhaus support}

\bigskip

\textbf{N. H. Bingham and A. J. Ostaszewski}

\bigskip

\textsl{In memoriam Herbert Heyer 1936-2018}

\bigskip
\end{center}

\noindent \textbf{Abstract. }Motivated by a Steinhaus-like interior-point
property involving the Cameron-Martin space of Gaussian measure theory, we
study a group-theoretic analogue, the Steinhaus triple $(H,G,\mu ),$ and
construct a Steinhaus support, a Cameron-Martin-like subset, $H(\mu )$ in
any Polish group $G$ corresponding to `sufficiently subcontinuous' measures $%
\mu ,$ in particular for `Solecki-type' reference measures.

\bigskip


\bigskip


\bigskip

\section{Introduction}

For many purposes, one needs a \textit{reference measure. }In discrete
situations such as the integers $\mathbb{Z}$, one has \textit{counting
measure. }In Euclidean space $\mathbb{R}^{d}$, one has \textit{Lebesgue
measure. }In locally compact groups, one has \textit{Haar measure. }In
infinite-dimensional settings such as Hilbert space, one has neither local
compactness nor Haar measure. Here various possibilities arise, some
pathological -- which we avoid by restricting to \textit{Radon} measures
(below). One is to use Christensen's concept of \textit{Haar-null sets }%
(below),\textit{\ }even though there is no Haar measure; see [Chr1,2],
Solecki [Sol], and the companion paper to this, [BinO7]. Another is to use 
\textit{Gaussian measures;} for background see e.g. Bogachev [Bog1,3], Kuo
[Kuo] and for Gaussian processes, Lifshits [Lif], Marcus and Rosen [MarR],
Ibragimov and Rozanov [IbrR]. Another is to weaken the invariance of measure
to (relative) \textit{quasi-invariance -- }see [Bog2, \S 9.11, p. 304-305],
and also \S 9.3 and \S 9.16 below.

Hilbert spaces are rather special, and the natural setting for Christensen's
Haar-null sets is Banach spaces. The Banach and Hilbert settings combine (or
intertwine) in Gross's concept of \textit{abstract Wiener space, }where
(identifying a Hilbert space $H$ with its dual) one has a triple $B\subseteq
H\subseteq B^{\ast \ast }$, with both inclusions continuous dense
embeddings. (It is tempting, but occasionally misleading, to speak of
embedded `subspaces' despite either subset here typically having a topology 
\textit{finer} than the \textit{subspace} topology induced by the containing
space, hence the use of quotation marks.) This is (essentially) the setting
of \textit{reproducing-kernel Hilbert spaces (RKHS); }see e.g. Berlinet and
Thomas-Agnon [BerTA]. Crucial here is the Cameron-Martin(-Maruyama-Girsanov)
theorem ([CamM1,2,3], [Gir]; [Bog1, 2.4], [Bog3, 1.4]). Here a suitable%
\textit{\ translation} gives a change of measure, the two measures being
equivalent, with Radon-Nikodym derivative given by the \textit{%
Cameron-Martin formula, }$(CM)$ below. Crucial also are the \textit{Gaussian
dichotomy} results (two Gaussian measures on the same space are either
equivalent or mutually singular). One has equivalence under translation%
\textit{\ }exactly when the translator is in the \textit{Cameron-Martin
space }[Bog1, 2.2]; these are the \textit{admissible translators}. We note
that the Girsanov \textit{change of measure }(by translation, using $(CM)$)
is the key to, e.g., Black-Scholes theory in mathematical finance (for
background see e.g. [BinK]).

Our purpose here is to construct a group-theoretic analogue of the
Cameron-Martin space arising in Gaussian measure theory. We are motivated by
a \textit{relativization of the Steinhaus} \textit{interior-point property}
[Ste], to be introduced below (important to classical regular variation --
see e.g. [BinGT, Th. 1.1.1]), i.e. with the notion of interior relativized
to a distinguished subset equipped with a finer topology:\ in brief a
`relativized version'. Though the classical paradigm may fail in an
infinite-dimensional Hilbert space, it can nevertheless hold relative to an
embedded (necessarily, compactly so -- see the final assertion of Th. 3.4)
space. Such is precisely the case when interiors are taken relative to the
Cameron-Martin space.

We recall the Gaussian context in a locally compact topological group $G$.
For simplicity, take $G$ Euclidean. Then matters split, according to the
support of the Gaussian measure. If this is the whole of $G$, the measure
has a density (given by the classical and familiar Edgeworth formula for the
multi-normal (multi-variate normal) distribution, see e.g. [BinF, 4.16]). If
not, the measure is singular viewed on $G$, but becomes non-singular when
restricted to the subgroup generated by its support (see \S 9.12). (This
situation is familiar in statistics: behaviour may seem degenerate only
because it is viewed in a context bigger than its natural one; see e.g.
[BinF].) Another instance of a similar `support-degeneracy' phenomenon
arises in the It\^{o}-Kawada theorem, when a (suitably non-degenerate)
probability measure $\mu $ has its convolution powers converging to Haar
measure on a subgroup $F$ of $G$, the closed subgroup generated by its
support [Hey1, \S\ 2.1]. In each case, the moral is the obvious one: if one
begins in the wrong context, identify the right one and start again.

Let $X$ be a locally convex topological vector space; it suffices for us to
take $X$ a separable\textit{\ Fr\'{e}chet }space (that is, having a
translation-invariant complete metric). Equip this with a Gaussian
(probability) measure $\gamma $ (`gamma for Gaussian', following [Bog1, Ch.
2] and \S 9.10; for Radon Gaussian measures in this context see [Bog1,
Ch.3]). Suppose further that the dual satisfies $X^{\ast }\subseteq
L^{2}(\gamma )$. Write $\gamma _{h}(K):=\gamma (K+h)$ for the translate by $%
h.$ \textit{Relative} \textit{quasi-invariance} of $\gamma _{h}$ and $\gamma
,$ that for all compact $K$ 
\[
\gamma _{h}(K)>0\text{ iff }\gamma (K)>0, 
\]%
holds relative to a set of vectors $h\in X$ (the \textit{admissible
translators}) forming a vector subspace known as the\textit{\ Cameron-Martin
space, }$H(\gamma )$. Then, in fact, $\gamma _{h}$ and $\gamma $ are
equivalent, $\gamma \sim \gamma _{h},$ iff $h\in H(\gamma ).$ Indeed, if $%
\gamma \sim \gamma _{h}$ fails, then the two measures are mutually singular, 
$\gamma _{h}\bot \gamma $ (the Hajek-Feldman Theorem -- cf. [Bog1, Th.
2.4.5, 2.7.2]).

Our key inspiration is that, for any non-null measurable subset $A$ of $X,$
the difference set $A-A$ contains a $|.|_{H}$-open nhd (neighbourhood) of $0$
in $H,$ i.e. $(A-A)\cap H$ contains a $H$-open nhd of $0$ -- see [Bog1, p.
64]. This flows from the continuity in $h$ of the density of $\gamma _{h}$
wrt $\gamma $ ([Bog1, Cor. 2.4.3]), as given in the \textit{%
Cameron-Martin-Girsanov formula}:%
\begin{equation}
\exp \left( \hat{h}(x)-\frac{1}{2}||\hat{h}||_{L^{2}(\gamma )}^{2}\right) 
\tag{$CM$}
\end{equation}%
(where $\hat{h}$ `Riesz-represents' $h,$ i.e. $x^{\ast }(h)=\langle x^{\ast
},\hat{h}\rangle ,$ for $x^{\ast }\in X^{\ast },$ as in \S 5). Thus here a
modified Steinhaus Theorem holds: the \textit{relative-interior-point theorem%
}.

In a locally compact topological group, Gaussian measures $\gamma $ may be
defined: see e.g. Heyer [Hey1, 5.2] (in the sense of Parthasarathy; cf.
[Hey1, 5.3] for Gaussianity in the sense of Bernstein). See also \S 9.12.
Such a $\gamma $ may be singular w.r.t. a (left) Haar measure $\eta $. Such
is the case in the Euclidean case, as above, with the Gaussian having its
support on a proper linear subspace $H$, and in this case $A-A$ with $A$
measurable and non-$\gamma $-null will only have non-empty \textit{relative} 
$H$-interior (and quasi-invariance only relative to $H$). We recall two
results due to Simmons [Sim] (cf. Mospan [Mos]; for generalizations beyond
the locally compact case, using results here, see also the companion paper
[BinO7]): (1) a measure $\mu $ is singular w.r.t. (left) Haar measure $\eta
_{G}$ on $G$ if and only if $\mu $ is concentrated on a $\sigma $-compact
subset $B$ such that $B^{-1}B$ has void interior (as in the Euclidean
example); (2) $\mu $ is absolutely continuous w.r.t. $\eta _{G}$ ($\mu \ll
\eta _{G})$ iff the group $G$ has the\textit{\ Steinhaus property}: for each
non-$\mu $-null compact set $K$,%
\[
1_{G}\in \mathrm{int}(K^{-1}K). 
\]%
(This does not preclude having $\mu (K)>0$ and $\mu (K+h)=0$ for some $K$
and $h.)$ Other characterizations of (Haar-) absolute continuity and
singularity are studied in [LiuR] and [LiuRW], and in the related
[Gow1,2,3]; see also [Pro] and [BarFF] on singularity. In certain locally
compact groups (e.g. [Hey1, 5.5.7] for the case $G=\mathbb{R}^{m}\times 
\mathbb{T}^{n}$ with $\mathbb{T}$ the unit circle) the condition $\mu \ll
\eta _{G}$ may imply that the support of a Gaussian probability measure $\mu 
$ is $G;$ see, however, [Hey1, 5.5.8] for an example of a Gaussian with full
support which is `Haar-singular'.

We develop an analogue of these relative interior results for a general
Polish group $G.$ This first leads, by analogy with an abstract Wiener space
triple [Bog1, 3.9], [Str, 4.2], to the concept of a \textit{Steinhaus triple}
$(H,G,\mu ),$ which we study in \S \S 2-4, demonstrating `relativized
variants' of classical results. In \S 5 we exhibit a link between the group
context and the classical Cameron-Martin theory above by verifying that a
divisible abelian group with an $\mathbb{N}$-homogeneous group-norm (below)
is in fact a topological vector space. In \S 6 we extend our usage in \S 3
of the notions of \textit{subcontinuity} and \textit{selective subcontinuity}
of a measure (introduced in [BinO4,7]), and in Theorem 6.1 establish the key
property of a \textit{Solecki reference measure}. Then in \S 7 for a given
Polish group $G$ and `sufficiently subcontinuous' measure $\mu $ we
construct a corresponding subset $H(\mu ),$ which together with $G$ and $\mu 
$ forms a Steinhaus triple (possibly `selective':\ see below). We verify
that it is an analogue of the Cameron-Martin space when $G$ is a Hilbert
space regarded as an additive group. In \S 8 we examine the extent of the
embedded `subspace' $H(\mu ).$ We close with complements in \S 9.

\section{Steinhaus triples: the context}

\qquad \textbf{Context.} Recall that a \textit{Polish} space $X$ is
separable and \textit{topologically complete}, i.e. its topology $\tau _{X}$
may be generated by a complete metric. Throughout the paper $G\ $will be a 
\textit{Polish group} [BecK]: a topological group which is a Polish space.
By the Birkhoff-Kakutani Theorem ([Bir], [Kak3]; cf. [DieS, \S 3.3]) we may
equip $G$ with a \textit{left-invariant} metric $d_{G}^{L}$ (equivalently,
with a (\textit{group-})\textit{norm} $||g||:=d_{G}^{L}(g,1_{G}),$ as in
[BinO2] -- `pre-norm' in [ArhT]) that generates its topology $\tau _{G}.$
(So $d_{G}^{L}(g,g^{\prime })=||g^{-1}g^{\prime }||$ and the corresponding
right-invariant metric is $d_{G}^{R}(g,g^{\prime }):=||g^{\prime }g^{-1}||$%
.) This metric, which is particularly useful, need not be complete, although 
$d=d_{G}^{L}+d_{G}^{R}$ is complete: see [TopH, Th. 2.3.5]. Nonetheless, the
group-norm endows $G$ with Fr\'{e}chet-like features, helpful here. When the
norm generating the topology is bi-invariant, \textit{Klee's completeness
theorem} [Kle], [DieS, Th. 8.16] asserts that if the topology is completely
metrizable, then in fact the norm is itself complete.

We fix a sequence of points $\{g_{n}\}_{n\in \mathbb{N}}$ dense in $G$; a
sequence $z_{n}\rightarrow 1_{G}$ will be called \textit{null}, and a null
sequence \textit{trivial} if it is ultimately constantly $1_{G}.$ For $%
\delta >0,$ by $B_{\delta }^{G}$ (resp. $B_{\delta }^{H}$) we denote the
open ball in $G$ (or $H)$ centered at $1_{G}$ of radius $\delta $ under $%
d_{G}^{L}$ (or $d_{H}^{L}$), which is \textit{symmetric}; by $\mathcal{B}(G)$
the \textit{Borel} sets; by $\mathcal{K}(G)$ the family of \textit{compact}
sets (carrying the Hausdorff metric induced by $d_{G}^{L}$).

Throughout, measure is to mean \textit{Borel measure} -- i.e. its domain
comprises the Borel sets of the relevant metrizable space [GarP] -- and such
a measure $\mu $ is \textit{Radon} if it is \textit{locally finite} (so that
each point has a neighbourhood of finite measure) and the relevant Borel
sets are \textit{inner compact regular}, i.e. ($\mu $-)approximable by
compacts from within:%
\[
\mu (B)=\sup \{\mu (K):K\in \mathcal{K}(G),K\subseteq B\}
\]%
([Sch] and \S 9.2); being locally finite on a separable metric space such a
measure is $\sigma $-finite. A $\sigma $-finite measure on a metric space is
necessarily \textit{outer regular} ([Kal, Lemma 1.34], cf. [Par, Th. II.1.2]
albeit for a probability measure), i.e. approximable by open sets from
without, and, when the metric space is completely metrizable, inner regular
([Bog2, II. Th. 7.1.7], [Xia, Th. 1.1.8], cf. [Par, Ths. II.3.1 and 3.2]).

For $\mu $ a Radon measure, we write $\mathcal{M}_{+}(\mu )$ for the
non-null sets, and put 
\[
\mathcal{K}_{+}(\mu ):=\mathcal{K}(G)\cap \mathcal{M}_{+}(\mu ). 
\]%
By $\mathcal{P}(G)$ we denote the family of (Borel) \textit{probability}
measures $\mu $ on $G,$ i.e. with $\mu (G)=1$ so these are Radon (as above);
by $\mathcal{U}(G)$ the \textit{universally measurable} sets (i.e.
measurable with respect to every measure $\mu \in \mathcal{P}(G)$ in the
sense of measure completion [Bog2.I, Th. 1.5.6, Prop. 1.5.11 `Lebesgue
completion'] -- for background and literature, see [BinO7]). We recall that $%
N\subseteq G$ is (left) \textit{Haar-null} in $G$ if, for some Borel $%
B\supseteq N$ and $\mu \in \mathcal{P}(G),$ $\mu (gB)=0$ for all $g\in G.$

\bigskip

\textbf{Definition. }Call $(H,G,\mu )$ a \textit{Steinhaus triple} if $G$ is
a Polish group with (group-) norm $||.||_{G}$ (notation as above), $\mu $ a
Radon measure on $G$ and $1_{G}\in H\subseteq G,$ a \textit{continuously
embedded subset} with group-norm $||.||_{H}$, having the property that for $%
K\in \mathcal{K}_{+}(\mu )$ there is in $H$ an ($H$-) open neighbourhood $U$
of $1_{G}$ (i.e. $U$ is open under $||.||_{H})$ such that%
\[
U\subseteq K^{-1}K. 
\]%
(So $\mathrm{cl}_{G}(U)$ is compact in $G.$) The latter condition links the
topological with the algebraic structure; the norm on $H$ introduces a
topology on $H$ \textit{finer} than the subspace topology induced by $G$ --
cf. [Bog1, Ch. 2], [BerTA], [Gro1,2], [Str, \S 4.2]. In view of its
distinguished status we will refer to $H$ as the \textit{Steinhaus support}.
Our aim is to establish structural similarities with the classical
Cameron-Martin space $H(\gamma )$ of \S 1 ([Bog1, \S 2.4], \S 3, \S 8 and \S %
9). Below, $\tau _{G}$ and $\tau _{H}$ will be the topologies of $G$ and $H;$
$\tau _{G|H}$ will be the topology induced on $H$ by $\tau _{G}.$ So 
\[
\tau _{G|H}\subseteq \tau _{H}. 
\]

It will be helpful to bear in mind the following illuminating example.

\bigskip

\textbf{Cautionary Example. }Consider the additive group $G=\mathbb{R}$ with
the Euclidean topology and its additive subgroup $H=\mathbb{Q}$ with \textit{%
discrete} topology. Enumerate $\mathbb{Q}$ as $\{q_{n}\}_{n\in \mathbb{N}}$
and put%
\[
\mu :=\tsum\nolimits_{n=1}^{\infty }2^{-n}\delta _{q_{n}}, 
\]%
with $\delta _{x}$ the Dirac measure at $x.$ Here $\tau _{G|H}\subseteq \tau
_{H}$. For $K\subseteq \mathbb{R}$ (Euclidean) compact, $\mu (K)>0$ iff $%
K\cap \mathbb{Q\neq \emptyset }$, and then $\{0\}$ ($\subseteq K-K$) is $%
\tau _{H}$-open.

\bigskip

The topological link between $H$ and $G$ above is at its neatest and most
thematic in norm language. But, as $H$ need not be a subgroup, it would
suffice for the continuous embedding to be determined by just a metric on $H$%
, or more generally a choice of refining topology. Furthermore, as in any
abstract definition of inner regularity, one is at liberty here to restrict
attention to a, possibly countable, subfamily of $\mathcal{K}_{+}(\mu )$ --
see \S 9.8. Such variants will be referred to below as \textit{selective
Steinhaus triples}. See the Remarks after Th. 7.1 and after Prop. 8.3 below.

\bigskip

\noindent \textbf{Remarks.} 1. The inclusion above implies that 
\[
K\cap Ku\neq \emptyset \qquad (u\in U) 
\]%
(indeed, if $u=k^{-1}k^{\prime }$, then $k^{\prime }=ku).$ In Theorem 3.4
below we strengthen this conclusion to yield the measure-theoretic
`Kemperman property', introduced in [Kem], as in [BinO4]. (One would expect
this to imply shift-compactness for $`G$-shifts of $H$', as indeed is so --
see Th. 3.8 below.)

2. Note that $W\cap H\ $is open in $H$ for $W$ open in $G$. It follows that
if $G$ is not locally compact, then, for $U\subseteq K^{-1}K$ as above, $U\ $%
is \textit{nowhere dense in} $G$ (as $\mathrm{cl}_{G}U$ is compact, its
interior in $G$ must be empty). Theorems 3.1 and 3.2 below expand this
remark to stark category/measure dichotomies.

\bigskip

We close this section with an observation which we need several times in the
next section.

\bigskip

\noindent \textbf{Lemma 2.1. }\textit{For }$(H,G,\mu )$\textit{\ a Steinhaus
triple, if} $H\ $\textit{is Polish in its own topology, then}%
\[
\mathcal{B}(H)\subseteq \mathcal{B}(G), 
\]%
\textit{\ i.e. if }$B$ \textit{is Borel in} $H,$ \textit{then }$B$ \textit{%
is Borel in }$G.$ \textit{In particular, if }$K\subseteq H$ \textit{is
compact in }$\tau _{H}$\textit{, then it is compact in }$\tau _{G}$\textit{: 
}$\mathcal{K}(H)\subseteq \mathcal{K}(G).$

\bigskip

\noindent \textbf{Proof. }As $H$ is Polish, $B$ being Borel is a Lusin space
(cf. [Sch, Ch. 2], [Rog, 1 \S\ 2.1]), so we may write $B$ as an injective
continuous image of the irrationals, $B=f(\mathbb{N}^{N})$ say, with $f:%
\mathbb{N}^{\mathbb{N}}\rightarrow H$ continuous (with $H$ under $\tau _{H}$%
) [Rog, 1. Cor 2.4.2]. Now the embedding $\iota :H\rightarrow G\ $is
continuous and so $B=\iota \circ f(\mathbb{N}^{\mathbb{N}})$ is an injective
continuous image of the irrationals, and so a Borel subset of $G$ [Rog, 1.
Th. 3.6.1]. The final assertion is clear (since $K=\iota (K)$ is compact, as 
$\tau _{G|K}\subseteq \tau _{H|K}$). $\square $

\section{Steinhaus triples: the general case}

In this section we study general topological properties of Steinhaus
triples, foremost among which is \textit{local quasi-invariance} (Theorem
3.5 below), a much weakened version of \textit{relative} quasi-invariance
(which we consider separately in the next section), i.e. relative to a
subgroup of `admissible' translators. This is preceded by a technical result
(Theorem 3.4) reminiscent of a lemma due to Kemperman [Kem] -- cf. [Kuc,
Lemma 3.7.2]; this will be revisited in another context in \S \S 6,7. We
close the section with Theorem 3.8, deducing a new property of the
Cameron-Martin space $H(\gamma ).$

For the sake of clarity, we emphasize that $\mathcal{K}_{+}(\mu )=\mathcal{K}%
(G)\cap \mathcal{M}_{+}(\mu )$: the compactness referred to is thus in the
sense of the topology on $G.$ Our first result confirms that the Steinhaus
support will always be meagre in our setting:

\bigskip

\noindent \textbf{Theorem 3.1.} \textit{For }$(H,G,\mu )$\textit{\ a
Steinhaus triple with }$H$\textit{\ Polish: if }$H$ \textit{is a dense
subgroup of }$G$\textit{, then either\newline
}\noindent i)\textit{\ }$H$\textit{\ is meagre in }$G$\textit{\ (so }$%
\mathrm{int}_{G}(H)=\emptyset ,$\textit{), or else\newline
}\noindent ii) $G$\textit{\ is locally compact and }$H=G$\textit{.}

\bigskip

\noindent \textbf{Proof.} Suppose (i) fails. Choose any $K\in \mathcal{K}%
_{+}(\mu )$ and consider any $\delta >0$ so small that $B_{\delta }^{H},$
the $H$-ball of radius $\delta >0$ centered at $1_{G},$ satisfies%
\[
B_{\delta }^{H}\cdot B_{\delta }^{H}=B_{2\delta }^{H}\subseteq K^{-1}K. 
\]%
For any countable dense subset $D$ of $H$%
\[
H=\bigcup\nolimits_{d\in D}d\cdot B_{\delta }^{H} 
\]%
(refer to the metric $d_{H}^{L}),$ so $d\cdot B_{\delta }^{H}$ is non-meagre
in $G$ for some $d\in D.$ As $G$ is a topological group, $B_{\delta
}^{H}=d^{-1}d\cdot B_{\delta }^{H}$ is also non-meagre. Now $B_{\delta
}^{H}, $ being open in $H,$ is Borel in $G$ by Lemma 2.1, and so has the
Baire property by Nikodym's Theorem [Rog, Part 1 \S 2.9], and is non-meagre
in $G.$ By the Piccard-Pettis Theorem ([Pic], [Pet], cf. [BinO2, Th. 6.5]),
applied in $G,$ there is $r>0$ with%
\[
B_{r}^{G}\subseteq B_{\delta }^{H}\cdot B_{\delta }^{H}=B_{2\delta
}^{H}\subseteq K^{-1}K. 
\]%
So $B_{2\delta }^{H}$ for all small enough $\delta >0$ contains $B_{r}^{G}$
for some $r>0,$ and $B_{r}^{G}$ has compact closure in $G.$ So $G=H$ (for
any $g\in G\ $choose $h\in gB_{r}^{G}\cap H;\ $then $g\in
hB_{r}^{G}\subseteq hB_{2\delta }^{H}\subseteq H$), the two topologies
coincide and $G$ is locally compact. $\square $

\bigskip

Next we turn from category to measure negligibility.

\bigskip

\noindent \textbf{Theorem 3.2.} \textit{For }$(H,G,\mu )$\textit{\ a
Steinhaus triple with }$H$\textit{\ Polish: if }$H\ $\textit{is a dense
subgroup of }$G$\textit{, then either\newline
}\noindent i)\textit{\ }$H$\textit{\ is }$\mu $\textit{-null, or else}%
\newline
\noindent ii)\textit{\ }$H$ \textit{is locally compact under its own
topology, with }$\mu _{H}\ll \eta _{H}$\textit{\ for }$\mu _{H}$\textit{\
the restriction of }$\mu $ \textit{to }$H$ \textit{and }$\eta _{H}$ \textit{%
a Haar-measure on }$H.$

\bigskip

\noindent \textbf{Proof.} Again suppose (i) fails and again recall that by
Lemma 2.1 the open subsets of $H$ have the Baire property as they are Borel
in $G$. Let $\mu _{H}$ denote the restriction of $\mu $ to the Borel subsets
of $H:$%
\[
\mu _{H}(B)=\mu (B\cap H)\qquad (B\in \mathcal{B}(H)). 
\]%
The resulting measure is still Radon: it is a locally finite (since $\tau
_{G|H}\subseteq \tau _{H}$) Borel measure and $H$ is Polish (cf. \S 2). So
since $\mu _{H}(H)=\mu (H)>0,$ there is a compact subset $K\subseteq H$
which is $\mu _{H}$-non-null. This set $K$ is compact also in the sense of $%
G $ (Lemma 2.1) and $\mu $-non-null. Hence again for some $\delta >0$%
\begin{equation}
B_{\delta }^{H}\subseteq K^{-1}K\subseteq H,  \tag{$KH$}
\end{equation}%
the latter inclusion as $H$ is a subgroup. So the topology of $H$ is locally
compact (and indeed $\sigma $-compact). Hence $H$ supports a left Haar
measure $\eta _{H}$ and so%
\[
\mu _{H}\ll \eta _{H}, 
\]%
by the Simmons-Mospan theorem ([Sim], [Mos]; [BinO7], \S 1), since $\mu _{H}$
has the Steinhaus property (\S 1) on $H.$ So 
\[
\mu (B)=\mu _{H}(B\cap H)+\mu (B\backslash H)\qquad (B\in \mathcal{B}(G)), 
\]%
with%
\[
\mu (B\cap H)=\int_{B\cap H}\frac{d\mu _{H}}{d\eta _{H}}d\eta _{H}(h).\qquad
\square 
\]

\bigskip

\noindent \textbf{Remarks.} 1. For $H$ locally compact, as above, the
subgroup $H$ of $G$ is capable of being generated by any non-null compact
subset $K$ of $H$. See $(KH)$ above; use a dense set of translates of $%
B_{\delta }^{H}.$

\noindent 2. For $L$ compact in $H,$ as $L\backslash B_{\delta }^{H}$ is
closed in $H$ it is also compact in $H,$ and so also in $G$ (via the
continuous embedding). Take $L:=\mathrm{cl}_{H}B_{\delta }^{H}\subseteq
K^{-1}K\subseteq H;$ then $L$ is a compact set of $H$ and so of $G.$ So the
set $W:=G\backslash (L\backslash B_{\delta }^{H})$ is open in $G,$ and so%
\[
B_{\delta }^{H}=L\backslash (L\backslash B_{\delta }^{H})=L\cap W. 
\]%
That is, the topology of $H$ on each subset of the form $L:=\mathrm{cl}%
_{H}B_{\delta }^{H}$ is induced by the topology of $G:$%
\[
\tau _{H|L}=\tau _{G|L}. 
\]%
So $H$ is a countable union of compact \textit{subspaces} of $G.$ Compare
the Cautionary Example above.

\noindent 3. Being locally compact, $H\ $is topologically complete and so is
an absolute $\mathcal{G}_{\delta }$. But that means only that it is a $%
\mathcal{G}_{\delta }$ subset in any space $X$ whereof it is a subspace:%
\[
\tau _{X|H}=\tau _{H}. 
\]

\bigskip

A complementary result follows, in which we need to assume that, for $\tau
_{G}$-compact subsets of $H,$ the mapping $m_{K}(h):=\mu (Kh)$ is continuous
on $H$ at $1_{G}$ relative to $\tau _{G|H},$ the topology induced by $G.$
(This is a relativized version of the global concept of \textit{mobility}
studied by A. van Rooij and his collaborators -- see e.g. [LiuR].) The proof
below is a relativized version of one in [Gow1]; we give it here as it is
short and thematic for our development.

\bigskip

\noindent \textbf{Theorem 3.2}$^{\prime }$ (cf. [Gow1]). \textit{For }$%
(H,G,\mu )$\textit{\ a Steinhaus triple with }$H$\textit{\ Polish and dense
in }$G,$ \textit{and }$G\ $\textit{not} \textit{locally compact}: \textit{if 
}$h\mapsto \mu (Kh)$\textit{\ is }$\tau _{G|H}$\textit{-continuous at }$%
1_{G} $ \textit{on }$H$\textit{\ for each }$K\in \mathcal{K}_{+}(\mu )$ 
\textit{lying in }$H,$ \textit{then }$\mu (H)=0.$

\bigskip

\noindent \textbf{Proof. }Suppose otherwise. Then, referring as in Th. 3.2
to the restriction $\mu _{H}(B)=\mu (B)$ for $B\in \mathcal{B}(H),$ there is 
$K\subseteq H\ $compact in $H\ $(so also compact in $G$) with $\mu (K)>0.$
Consider any $\delta >0$ with $0<\delta <\mu (K).$ Put $\varepsilon :=\delta
/3.$ Choose open $U$ in $G$ with $K\subseteq U$ such that $\mu (U)<\mu
(K)+\varepsilon ,$ and $\eta >0$ so that%
\[
|\mu (Kh)-\mu (K)|\leq \varepsilon 
\]%
for $h\in B_{\eta }^{G}\cap H.$ W.l.o.g. $KB_{\eta }^{G}\subseteq U,$ whence 
$\mu (Kh\backslash K)\leq \mu (U\backslash K)\leq \varepsilon ,$ for $h\in
B_{\eta }^{G}\cap H,$ and likewise $\mu (K\backslash Kh)\leq 2\varepsilon ,$
because%
\[
\mu (K)-\varepsilon \leq \mu (Kh)\leq \mu (Kh\cap K)+\mu (U\backslash K). 
\]%
Combining yields%
\[
\mu (Kh\triangle K)\leq 3\varepsilon =\delta . 
\]%
So, since $\delta <\mu (K)$, 
\[
B_{\eta }^{G}\cap H\subseteq \{h\in H:\mu (Kh\triangle K)<\delta \}\subseteq
\{h\in H:Kh\cap K\neq \emptyset \}\subseteq K^{-1}K: 
\]%
\[
H\cap B_{\eta }^{G}\subseteq K^{-1}K. 
\]%
So%
\[
\mathrm{cl}_{G}(B_{\delta }^{G})=\mathrm{cl}_{G}(H\cap B_{\delta
}^{G})\subseteq K^{-1}K. 
\]%
Hence $G$ is locally compact, a contradiction. $\square $

\bigskip

\noindent \textbf{Theorem 3.3. }\textit{For }$(H,G,\mu )$\textit{\ a
Steinhaus triple, with }$H$\textit{\ Polish, if }$H$ \textit{is a dense
proper subgroup of }$G$\textit{, then }$H$\textit{\ is generically left
Haar-null in }$G$\textit{\ -- left Haar-null for quasi all }$\mu ^{\prime
}\in \mathcal{P}(G),$\textit{\ in the sense of the L\'{e}vy metric on }$%
\mathcal{P}(G).$\textit{\ In particular, this is so for the Cameron-Martin
space }$H(\gamma ).$

\bigskip

\noindent \textbf{Proof.} This follows from a result of Dodos [Dod, Cor. 9]
since $H$ is an analytic subgroup (cf. Lemma 2.1) with empty interior (Th.
3.1) and has $\mu $-measure zero (Th. 3.2$^{\prime }$). $\square $

\bigskip

We include here a similar result which is thematic.

\bigskip

\noindent \textbf{Theorem 3.3}$^{\prime }$\textbf{\ (Smallness). }\textit{%
For }$X$\textit{\ a Fr\'{e}chet space carrying} \textit{a Radon Gaussian
measure }$\gamma $\textit{,}\textbf{\ }$H(\gamma )$\textit{\ is
(generically) left Haar-null -- left Haar-null for quasi all }$\mu \in 
\mathcal{P}(X),$\textit{\ in the sense of the L\'{e}vy metric on }$\mathcal{P%
}(X).$

\bigskip

\noindent \textbf{Proof. }For $\gamma $ a Radon probability measure, $%
H:=H(\gamma )$ is a separable Hilbert space ([Bog1, 3.2.7], cf. [Bog1, p.
62]). So it is complete under its own norm, so a Polish space; it is
continuously embedded in $X$ [Bog1, 2.4.6] so, as a subset of $X,$ it is
analytic as in Lemma 2.1 (here:$\ \tau _{X|H}\subseteq \tau _{H}$). As $H$
has empty interior in $X$ (Th. 3.1), again by the result of Dodos [Dod, Cor.
9], $H$ is generically left-Haar null. $\square $

\bigskip

\noindent \textbf{Remark. }In the setting above, $H:=H(\gamma )$ is in fact
a $\sigma $-\textit{compact} subset of $X$: by [Bog1, 2.4.6], the $H$-closed 
\textit{unit ball} $U_{H}$ is weakly closed in $X\ $and, being convex (by
virtue of its norm), it is closed in $X$ [Rud2, 3.12], cf. [Con, \S 5.12];
but, as in the next theorem (Th. 3.4), it is a subset of some compact set of
the form $K-K$. So $U_{H}$ is itself compact in $X$ -- see [Bog1, 3.2.4].

Of course, if $X\ $is an infinite-dimensional space Hilbert space and $H$ is
a $\sigma $-compact subset of $X$, then, by Baire's theorem, $H\ $must have
empty interior in $X.$

\bigskip

\noindent \textbf{Theorem 3.4.} \textit{For }$(H,G,\mu )$\textit{\ a
Steinhaus triple and} $K\in \mathcal{K}_{+}(\mu ),$ \textit{there are }$%
\varepsilon ,\delta >0$\textit{\ such that}%
\[
\mu (K\cap Kh)\geq \delta \text{ for all }h\in H\ \text{with }||h||_{H}\leq
\varepsilon . 
\]%
\textit{In particular,}%
\[
\mu _{-}^{H}(K):=\sup_{\varepsilon >0}\inf \{\mu (Kh):h\in B_{\varepsilon
}^{H}\}\geq \delta , 
\]%
\textit{so that for any} \textit{null sequence} $\mathbf{t}=\{t_{n}\}$ 
\textit{in }$H$ \textit{(i.e. with} $t_{n}\rightarrow _{H}1_{G}),$ 
\[
\mu _{-}^{\mathbf{t}}(K):=\lim \inf\nolimits_{n\rightarrow \infty }\mu
(Kt_{n})\geq \delta ; 
\]%
\textit{furthermore, for some} $r>0,$%
\[
K\cap Kh\in \mathcal{M}_{+}(\mu )\qquad (h\in B_{r}^{H}):\qquad
B_{r}^{H}\subseteq K^{-1}K. 
\]

\bigskip

\noindent \textbf{Proof. }Suppose otherwise. Then for some $K\in \mathcal{K}%
_{+}(\mu )$ and for each pair $\varepsilon ,\delta >0$ there is $h\in H$
with $||h||_{H}<\varepsilon $ and%
\[
\mu (K\cap Kh)<\delta . 
\]%
So in $H$ there is a sequence $t_{n}\rightarrow _{H}1_{G}$ with 
\[
\mu (K\cap Kt_{n})<2^{-n-1}\mu (K) 
\]%
for each $n\in \mathbb{N}.$ Take%
\[
M:=K\cap \dbigcup\nolimits_{n\in \mathbb{N}}Kt_{n}. 
\]%
Then%
\[
\mu (M)<\mu (K)/2. 
\]%
So we may choose a compact $\mu $-non-null $K_{0}\subseteq K\backslash M;$
then, since $(H,G,\mu )$ is a Steinhaus triple, there is in $H$ a non-empty
open nhd $V$ of $1_{G}$ with 
\[
V\subseteq K_{0}^{-1}K_{0}. 
\]%
Now $t_{m}\in V$ for all large $m$. Fix such an $m;$ then%
\[
K_{0}\cap K_{0}t_{m}\neq \emptyset . 
\]%
Consider $k_{0},k_{1}\in K_{0}$ with $k_{0}=k_{1}t_{m}\in K_{0}\cap
K_{0}t_{m}\subseteq K\cap \dbigcup\nolimits_{n}Kt_{n}=M;$ as $K_{0}$ is
disjoint from $M,$ this is a contradiction. $\square $

\bigskip

An immediate corollary is

\bigskip

\noindent \textbf{Theorem 3.5 (Local quasi-invariance).} \textit{For }$%
(H,G,\mu )$ \textit{a Steinhaus triple and }$B\in \mathcal{B}(G)$\textit{:
if }$\mu (B)>0,$\textit{\ then there exists }$\delta >0$ \textit{so that }$%
\mu (Bh)>0$\textit{\ for all }$h\in H$\textit{\ with }$||h||_{H}\leq \delta
. $

\bigskip

\noindent \textbf{Proof. } This follows from Th. 3.4 since there is compact $%
K\subseteq B$ with $\mu (K)>0.$ Then for all sufficiently small $h\in H$ $%
\mu (Bh)\geq \mu (Kh)>0.$ $\square $

\bigskip

\noindent \textbf{Corollary 3.1.} \textit{For each null sequence }$%
t_{n}\rightarrow 1_{G}$\textit{\ in }$H$\textit{\ and }$\mu $\textit{%
-non-null }$K,$%
\[
0<\mu _{-}(K)\leq \lim \inf \mu (t_{n}K)\leq \lim \sup \mu (t_{n}K)\leq \mu
(K), 
\]%
\textit{and so, for }$\mu $\textit{-non-null compact }$L\subseteq K,$%
\[
0<\mu _{-}(L)\leq \mu _{-}(K)\text{ and }\mu (L)\leq \mu (K). 
\]

\bigskip

\noindent \textbf{Proof.} Writing $\mu ^{\delta }(K):=\inf \{\mu (Kh):h\in
B_{\delta }^{H}\},$ for any $\delta >0$ and all large enough $n$%
\[
\mu ^{\delta }(K)\leq \mu (t_{n}K)\leq \mu ^{\delta }(K)+\delta , 
\]%
yielding the lower bound when $\delta \downarrow 0$. Also for $\delta >0,$
there is open $U\supseteq K$ with $\mu (U\backslash K)\leq \delta ,$ and so
as $K$ is compact for all large enough $n$%
\[
\mu (t_{n}K)\leq \mu (K)+\delta . 
\]

For each $\delta >0$ choose $t_{\delta }\in B_{\delta }^{H}$ with $\mu
^{\delta }(K)\leq \mu (t_{\delta }K)\leq \mu ^{\delta }(K)+\delta .$ Then,
for non-null $L\subseteq K,$ since $t_{\delta }L\subseteq t_{\delta }K$, by
the earlier proved assertions,%
\[
\mu _{-}(L)\leq \lim \inf \mu (t_{\delta }L)\leq \lim \inf \mu (t_{\delta
}K)\leq \lim \inf [\mu ^{\delta }(K)+\delta ]=\mu _{\_}(K), 
\]%
and%
\[
\lim \sup \mu (t_{\delta }L)\leq \mu (L)\leq \mu (K).\qquad \square 
\]

\noindent \textbf{Remarks.} 1. Above, if $\mu (t_{n}K)\rightarrow \mu
_{0}\geq \mu _{-}(K),$ then%
\[
\mu (\tbigcap\nolimits_{n\in \mathbb{N}}(K\backslash t_{n}K))=\mu (K)-\mu
_{0}\leq \mu (K)-\mu _{-}(K). 
\]

\noindent 2. Evidently, for disjoint non-null compact $K,L$ 
\[
\mu (K)+\mu (L)=\mu (K\cup L)\geq \mu _{-}(K\cup L)\geq \mu _{-}(K)+\mu
_{-}(L). 
\]%
If one of these is sharp, one may imagine passsing through a subsequence $%
K_{n}$ with $\mu (K_{n})>\mu _{-}(K_{n})>0.$

\bigskip

There is no Steinhaus-like assumption on the measure $\mu $ in the following
result, which, standing in apposition to Th. 3.4, is a kind of converse.

\bigskip

\noindent \textbf{Proposition 3.1 }([BinO7, L. 1]). \textit{For }$H\subseteq
G$ \textit{continuously embedded in} $G$, $\mu \in \mathcal{P}(G)$\textit{\
and }$K\in \mathcal{K}_{+}(\mu )$: \textit{if }$\mu _{-}^{H}(K)>0$\textit{,\
then there is }$\delta >0$\textit{\ with }%
\[
\Delta /4\leq \mu (K\cap Kt)\qquad (t\in B_{\delta }^{H}), 
\]%
\textit{so that}%
\begin{equation}
K\cap Kt\in \mathcal{M}_{+}(\mu )\qquad (t\in B_{\delta }^{H}). 
\tag{$\ast
$}
\end{equation}%
\textit{In particular,}%
\[
K\cap Kt\neq \emptyset \qquad (t\in B_{\delta }^{H}), 
\]%
\textit{or, equivalently,}%
\begin{equation}
B_{\delta }^{H}\subseteq K^{-1}K,  \tag{$\ast \ast $}
\end{equation}%
\textit{so that }$B_{\delta }^{H}$ \textit{has compact closure under }$\tau
_{G}$\textit{.}

\bigskip

\noindent \textbf{Proof.} Put $H_{t}:=K\cap Kt\subseteq K$. Take $\Delta
:=\mu _{-}^{H}(K)>0.$ Then for any small enough $\delta >0$, $\mu
(Kt)>\Delta /2$ for $t\in B_{\delta }^{H}.$ Fix such a $\delta >0.$

By outer regularity of $\mu $, choose $U$ open with $K\subseteq U$ and $\mu
(U)<\mu (K)+\Delta /4.$ By upper semicontinuity of $t\mapsto Kt$, w.l.o.g. $%
KB_{\delta }\subseteq U$ $.$ For $t\in B_{\delta }^{H},$ by finite
additivity of $\mu ,$ since $\Delta /2<\mu (Kt)$ 
\begin{eqnarray*}
\Delta /2+\mu (K)-\mu (H_{t}) &\leq &\mu (K)+\mu (Kt)-\mu (H_{t})=\mu (K\cup
Kt) \\
&\leq &\mu (U)\leq \mu (K)+\Delta /4.
\end{eqnarray*}%
Comparing extreme ends of this chain of inequalities gives%
\[
0<\Delta /4\leq \mu (H_{t})\qquad (t\in B_{\delta }^{H}). 
\]

For $t\in B_{\delta }^{H},$ as $K\cap Kt\in \mathcal{K}_{+}(\mu )$, take $%
s\in K\cap Kt\neq \emptyset ;$ then $s=at$ for some $a\in K,$ so $%
t=a^{-1}s\in K^{-1}K.$ Conversely, $t\in B_{\delta }^{H}\subseteq K^{-1}K$
yields $t=a^{-1}a^{\prime }$ for some $a,a^{\prime }\in K;$ then $a^{\prime
}=at\in K\cap Kt$. $\square $

\bigskip

L\"{o}wner showed in 1939 that there exists no ($\sigma $-finite)
translation-invariant measure on an infinite-dimensional Hilbert space
([Loe, \S 1], [Neu2]). This is contained in the result below: with $G$ a
Hilbert space regarded as an additive group and $\mu $ Radon, if $\mu $ is
translation-invariant, then $\mu (K)=\mu _{-}^{G}(K)>0$ for some compact $K,$
so $G$ is locally compact and so finite-dimensional.

\bigskip

\noindent \textbf{Corollary 3.2 }(cf. [Gow1])\textbf{.} \textit{If }$\mu
_{-}^{G}(K)>0$\textit{\ for some }$K\in \mathcal{K}_{+}(\mu ),$\textit{\
then }$G\ $\textit{is locally compact.}

\bigskip

\noindent \textbf{Proof.} If $\mu _{-}(K)>0,$ then there is $\delta >0$ such
that $\mu (tK)>\mu _{-}(K)/2$ for all $t\in B_{\delta }=B_{\delta }^{G};$ so
for $\Delta :=\mu _{-}(K)/2$%
\[
B_{\delta }\subseteq B_{\delta }^{\Delta }=\{z\in B_{\delta }:\mu
(Kz)>\Delta \}\subseteq K^{-1}K, 
\]%
and so $B_{\delta }$ has compact closure. $\square $

\bigskip

\noindent \textbf{Theorem 3.6.} \textit{For }$(H,G,\mu )$\textit{\ a
Steinhaus triple and} $K\in \mathcal{K}_{+}(\mu ),$ \textit{the set}%
\[
\mathcal{O}(K):=\{h\in H:\mu (Kh)>0\} 
\]%
\textit{is }$H$\textit{-open, and so }$\mu $\textit{\ is continuous on a
dense }$\mathcal{G}_{\delta }$\textit{\ of }$\mathcal{O}(K)$ \textit{(so off
a meagre set), i.e.}%
\[
\mu (Kh)=\mu _{-}^{H}(Kh)>0\qquad (\text{quasi all }h\in \mathcal{O}(K)). 
\]%
\textit{Conversely, for }$H,G\ $\textit{topological groups with }$H$\textit{%
\ a continuously embedded subgroup of }$G\ $\textit{and }$\mu \in \mathcal{P}%
(G)$\textit{: if }$\mathcal{O}(K)$ \textit{is open for }$K\in \mathcal{K}%
_{+}(\mu )$ \textit{and the above relative continuity property of }$\mu $%
\textit{\ holds on a dense} $\mathcal{G}_{\delta }$ \textit{subset of} $%
\mathcal{O}(K)$\textit{\ in }$H,$ \textit{then }$(H,G,\mu )$\textit{\ is a
Steinhaus triple. }

\bigskip

\noindent \textbf{Proof.} The first assertion follows from Th. 3.5 on local
quasi-invariance applied to $Kh$ with $h\in \mathcal{O}(K).$ Since the map $%
g\mapsto \mu (Kg)$ is upper semi-continuous (see. e.g. [BinO7, Prop. 1],
[Hey1, 1.2.8]), the second assertion follows from the first by the theorem
of Fort [For, R1] (cf. [Xia, Appendix I, Lemma I.2.2]) that an upper
semi-continuous map is continuous on a $H$-dense $\mathcal{G}_{\delta }$ in $%
\mathcal{O}(K).$

As for the converse, for $K\in \mathcal{K}_{+}(\mu ),$ since $1_{G}\in 
\mathcal{O}(K)$ there is $h\in H$ with $\mu (Kh)=\mu _{-}^{H}(Kh)>0.$ It now
follows by Th. 3.4 that $B_{r}^{H}\subseteq (Kh)^{-1}Kh$ for some $r>0,$ and
so $1_{G}\in hB_{r}^{H}h^{-1}\subseteq K^{-1}K,$ i.e. $1_{G}$ is an interior
point of $K^{-1}K$ (since $hB_{r}^{H}h^{-1}$ is open, as $H$ is a
topological group). $\square $

\bigskip

As a corollary of Th. 3.4, we now obtain a result concerning embeddability
into non-negligible sets (here the non-null measurable sets) of some
translated subsequence of a given null sequence. This property, first used
implicitly by Banach [Ban1,2], has been studied in various general contexts
by many authors, most recently under the term `shift-compactness' -- see
e.g. [Ost1]. The new context of a Steinhaus triple is notable in limiting
the null sequences to the distinguished subspace. Here the statement calls
for the passage from a null sequence in $H$ to its inverse sequence; this
inversion is of course unnecessary if $H^{-1}=H,$ e.g. if $H$ is a subgroup
of $G,$ as will be the case in Theorem 3.8 below. (The group-theoretic
approach to shift-compactness is that of a group action, here of translation
in $G$ -- see [MilO]; for applications see [Ost2].)

\bigskip

\noindent \textbf{Theorem 3.7 (Shift-compactness Theorem for Steinhaus
triples).} \textit{For }$(H,G,\mu )$\textit{\ a Steinhaus triple, }$\mathbf{h%
}$ \textit{a null sequence in }$H,$ \textit{and }$E\in \mathcal{M}_{+}(\mu
): $\textit{\ for }$\mu $\textit{-almost all }$s\in E$ \textit{there exists
an infinite }$\mathbb{M}_{s}\subseteq \mathbb{N}$ \textit{with }%
\[
\{sh_{m}^{-1}:m\in \mathbb{M}_{s}\}\subseteq E. 
\]

\bigskip

\noindent \textbf{Proof.} Fix a compact $K_{0}\subseteq E$ with $\mu
(K_{0})>0.$ Choose inductively a sequence $m(n)\in \mathbb{N}$ and
decreasing compact sets $K_{n}\subseteq K_{0}\subseteq E$ with $\mu
(K_{n})>0 $ such that%
\[
\mu (K_{n}\cap K_{n}h_{m(n)})>0. 
\]%
To check the inductive step, suppose $K_{n}$ already defined. As $\mu
(K_{n})>0,$ by Th. 3.4 there are $\delta ,\varepsilon >0$ such that 
\[
\mu (K_{n}\cap K_{n}h)\geq \delta \text{ for all }h\in H\ \text{with }%
||h||_{H}\leq \varepsilon . 
\]%
So there is $m(n)>n$ with $\mu (K_{n}\cap K_{n}h_{m(n)})>0.$ Putting $%
K_{n+1}:=K_{n}\cap K_{n}h_{m(n)}\subseteq K_{n}$ completes the inductive
step, and so the induction.

By compactness, select $s$ with%
\[
s\in \bigcap\nolimits_{m\in \mathbb{N}}K_{m}\subseteq K_{n+1}=K_{n}\cap
K_{n}h_{m(n)}\qquad (n\in \mathbb{N}). 
\]%
Choosing $k_{n}\in K_{n}\subseteq K$ with $s:=k_{n}h_{m(n)}$ gives $s\in
K_{0}\subseteq E,$ and%
\[
sh_{m(n)}^{-1}=k_{n}\in K_{n}\subseteq K_{0}\subseteq E. 
\]%
Finally take $\mathbb{M}:=\{m(n):n\in \mathbb{N}\}.$

As for the final assertion, recalling from \S 2 that $\mathcal{U}(G)$
denotes the universally measurable sets, define 
\[
F(H):=\bigcap\nolimits_{n\in \mathbb{N}}\bigcup\nolimits_{m>n}H\cap
Hh_{m}\qquad (H\in \mathcal{U}(G)). 
\]%
Then $F:\mathcal{U}(G)\rightarrow \mathcal{U}(G)$ and $F$ is monotone $%
(F(S)\subseteq F(T)\ $for $S\subseteq T);$ moreover, $s\in F(H)$ iff $s\in H$
and $sh_{m}^{-1}\in H$ for infinitely many $m$. It suffices to show that $%
E_{0}:=E\backslash F(E)$ is $\mu $-null (cf. the Generic Completeness
Principle [BinO1, Th. 3.4]). Suppose otherwise. Then, as $\mu (E_{0})>0,$
there exists a compact $K_{0}\subseteq E_{0}$ with $\mu (K_{0})>0.$ But
then, as in the construction above, $\emptyset \neq F(K_{0})\cap
K_{0}\subseteq F(E)\cap E_{0},$ contradicting $F(E)\cap E_{0}=\emptyset .$ $%
\square $

\bigskip

\noindent \textbf{Corollary 3.3.} \textit{If the subsequence embedding
property of Theorem 3.7 holds for all the null sequences in a set }$H$%
\textit{\ which is continuously embedded in }$G$\textit{\ for all }$E\in 
\mathcal{K}_{+}(\mu )$, \textit{then }$(H,G,\mu )$\textit{\ is a Steinhaus
triple. In particular, for any }$E\in \mathcal{M}_{+}(\mu )$ \textit{the set 
}$E^{-1}E$ \textit{has non-empty }$H$\textit{-interior.}

\bigskip

\noindent \textbf{Proof.} If in $H$ there is no open subset $U$ with $%
U\subseteq K^{-1}K,$ then there exists $h_{n}\in H\ $with $h_{n}\notin
B_{1/n}^{H}\backslash (K^{-1}K).$ Then there is $s\in K$ and an infinite $%
\mathbb{M}_{s}\subseteq \mathbb{N}$ with%
\[
\{sh_{m}^{-1}:m\in \mathbb{M}_{s}\}\subseteq K. 
\]%
So for any $m\in \mathbb{M}_{s},$ $h_{m}s^{-1}\in K^{-1}$, i.e. $h_{m}\in
K^{-1}K,$ a contradiction. $\square $

\bigskip

Another immediate corollary is the following result, which was actually our
point of departure.

\bigskip

\noindent \textbf{Theorem 3.8 (Shift-compactness Theorem for the
Cameron-Martin Space). }\textit{For }$X$\textit{\ a Fr\'{e}chet space
carrying a Radon Gaussian measure }$\gamma $\textit{\ with }$X^{\ast
}\subseteq L^{2}(\gamma )$\textit{, and }$H(\gamma )$\textit{\ the
Cameron-Martin space: if} $\mathbf{h}$ \textit{is null in }$H(\gamma ),$ 
\textit{and }$E\in \mathcal{M}_{+}(\gamma )$\textit{, then for }$\gamma $%
\textit{-almost all }$s\in E$ \textit{there exists an infinite }$\mathbb{M}%
_{s}\subseteq \mathbb{N}$ \textit{with }%
\[
\{s+h_{m}:m\in \mathbb{M}_{s}\}\subseteq E. 
\]

\bigskip

\noindent \textbf{Proof.} Regarding $X$ as an additive group, $(H(\gamma
),X,\gamma )$ is a Steinhaus triple, by [Bog1, p. 64]. As $H(\gamma )$ is a
subspace of $X$, $(-h_{n})$ is also a null sequence in $H(\gamma );$ by
Theorem 3.7, for $\gamma $-almost all $s\in E$ there is $\mathbb{M}%
_{s}\subseteq \mathbb{N}$ with 
\[
\{s-(-h_{m}):m\in \mathbb{M}_{s}\}\subseteq E.\qquad \square 
\]

\section{Steinhaus triples: the quasi-invariant case}

We begin with the definition of measure `relative quasi-invariance' promised
in \S 3. Classical results on this topic are given in the setting of
topological vector spaces -- see e.g. [GikS, Ch. VII], [Sko], [Yam2] -- with
the exception of [Xia] who develops the associated harmonic analysis in its
group setting. We adopt a similar approach here in order to pursue some
parallels with Cameron-Martin theory.

\subsection{Relative quasi-invariance}

\noindent \textbf{Definition. }Say that $\mu $ is \textit{relatively
quasi-invariant w.r.t }$H$, or just $H$\textit{-quasi-invariant} if $\mu
(hB)=0$ for all $\mu $-null Borel $B\in \mathcal{B}(G)$ (equivalently, $\mu $%
-null compact $B)$ and $h\in H$ $.$

\bigskip

Recall that $\mathrm{supp}_{G}(\mu )$ denotes the \textit{topological support%
}, which is the smallest closed set of full $\mu $-measure; for $\mu $
Radon, such a smallest closed set is guaranteed to exist [Bog2, Prop. 7.2.9].

\bigskip

\noindent \textbf{Proposition 4.1 }(cf. [Bog1, 3.6.1])\textit{. For }$%
(H,G,\mu )$ \textit{a Steinhaus triple with }$H$\textit{\ a subgroup and }$%
\mu \in \mathcal{P}(G),$\textit{\ a }$H$\textit{-quasi-invariant (Radon)
measure: if }$\mu (\mathrm{cl}_{G}H)=1,$ \textit{then\ }%
\[
S_{\mu }:=\mathrm{supp}_{G}(\mu )=\mathrm{cl}_{G}H.\mathit{\ } 
\]%
\textit{In particular, this is so for }$H$\textit{\ the Cameron-Martin space 
}$H(\gamma ).$

\bigskip

\noindent \textbf{Proof. }Let $L:=\mathrm{cl}_{G}H;$ then $L\subseteq S_{\mu
}.$ If the inclusion were proper: take $x\in L\backslash S_{\mu }.$ There is 
$V$ open in $G$ with $x\in V$ and $\mu (V)=0.$ But, as $x\in L,$ there is $%
h_{0}\in V\cap H,$ and so $1_{G}\in W:=h_{0}^{-1}V$ with $\mu (W)=0.$
W.l.o.g. $W^{-1}=W$ (otherwise pass to $W\cap W^{-1},$ which contains $%
1_{G}).$ So also $\mu (hW)=0$ for $h\in H$ (by $H$-quasi-invariance). Then,
by the definition of the support,%
\[
HW=\bigcup\nolimits_{h\in H}(hW)\subseteq X\backslash S_{\mu }, 
\]%
and so $\mu (HW)=0.$ But $\mathrm{cl}_{G}H\subseteq HW,$ for if the point $%
x\in \mathrm{cl}_{G}H,$ then its nhd $xW$ meets $H,$ in $h$ say; then $x\in
hW^{-1}=hW$ (since $xw=h$ implies $x=hw^{-1}$). So $\mu (\mathrm{cl}%
_{G}H)=0, $ contradicting the fact that $\mu (\mathrm{cl}_{G}H)=1$. $\square 
$

\bigskip

\noindent \textbf{Definition. }Following [Bog1, 3.6.2], say that for a
Steinhaus triple $(H,G,\mu )$ the measure $\mu $ is \textit{non-degenerate}
iff $S_{\mu }=\mathrm{cl}_{G}H.$

\bigskip

We close this subsection by tracing a measure-to-category dependence.

\bigskip

\noindent \textbf{Proposition 4.2 (From measure to category: nullity to
empty interior).}\textit{\ For }$(H,G,\mu )$ \textit{a Steinhaus triple with 
}$H$ \textit{a subgroup and }$\mu \in \mathcal{P}(G)$\textit{\ an }$H$%
\textit{-quasi-invariant} \textit{non-degenerate (Radon) measure: if }$\mu
(H)=0$\textit{, then}%
\[
\mathrm{int}_{G}(H)=\emptyset . 
\]%
\textit{In particular, this is so for }$H$ \textit{the Cameron-Martin space }%
$H(\gamma ).$

\bigskip

\noindent \textbf{Proof. }Suppose not. Then w.l.o.g. $1_{G}\in W:=$\textrm{%
int}$_{G}(H)\subseteq H$ (as $H$ is a subgroup and so $Ww^{-1}\subseteq H$
for $w\in W$). Also w.lo.g. $W=W^{-1}$ (otherwise pass to the non-empty open
subset $W^{-1}\cap W\subseteq H).$ But $\mu (H)=0,$ so $\mu (W)=0,$ and so
also $\mu (hW)=0$ for $h\in H$ (by $H$-quasi-invariance). Then, by the
definition of the support,%
\[
HW=\bigcup\nolimits_{h\in H}hW\subseteq X\backslash S_{\mu }. 
\]%
So $\mu (HW)=0,$ and the rest of the proof is as in Prop. 4.1 using $S_{\mu
}=\mathrm{cl}_{G}H$. $\square $

\subsection{Admissible translators for $\protect\mu $-quasi-invariance}

We close with a study of the algebraic structure of admissible translators
by considering the natural candidates for the Steinhaus support of a (Radon)
measure and a corresponding natural complement (inspired by the
Hajek-Feldman Dichotomy Theorem -- cf. [Bog1, Th. 2.4.5, 2.7.2]). The
results here (in particular Prop.4.5) will be used in \S 8. In what follows
the use of $Q$ (`q for quasi-invariance') is justified in Prop. 4.4 below;
write $\mu _{g}(B):=\mu (Bg)$ for $B\in \mathcal{B}(G),$ and put 
\begin{eqnarray*}
Q_{R} &=&Q_{R}(\mu ):=\{g\in G:(\forall K\in \mathcal{K}_{+}(\mu ))\text{ }%
\mu (Kg)>0\}, \\
Q &=&Q(\mu ):=\{g\in G:(\forall K\in \mathcal{K}_{+}(\mu ))[\mu (Kg)>0\text{
\& }\mu (Kg^{-1})>0]\}, \\
\mathcal{N} &\mathcal{=}&\mathcal{N}(\mu ):=\{g\in G:\mu _{g}\bot \mu
\},\qquad G_{0}=G_{0}(D):=\dbigcup\nolimits_{d\in D}d\mathcal{N}(\mu ),
\end{eqnarray*}%
where $D:=\{g_{n}:n\in \mathbb{N}\}$ is a dense subset of $G.$ Evidently 
\textbf{\ }%
\[
\mathcal{N}(\mu )=\dbigcap\nolimits_{K\in \mathcal{K}_{+}(\mu )}\{g:\mu
(Kg)=0\}=\dbigcap\nolimits_{n\in \mathbb{N}}\dbigcap\nolimits_{K\in \mathcal{%
K}_{+}(\mu )}\{g:\mu (Kg)<1/n\}. 
\]%
The definitions imply that 
\[
Q\subseteq Q_{R}\subseteq G\backslash \mathcal{N}(\mu ). 
\]

\bigskip

\noindent \textbf{Lemma 4.1.} $Q$ \textit{is a subgroup of} $G.$

\bigskip

\noindent \textbf{Proof. }First, $g\in Q$ iff $g^{-1}\in Q.$ Next take $%
x,y\in Q$ and $K\in \mathcal{K}_{+}(\mu ).$ As $Kx\in \mathcal{K}_{+}(\mu )$
and $y\in Q,$ $\mu (Kxy)=\mu ((Kx)y)>0.$ $\square $

\bigskip

\noindent \textbf{Proposition 4.3. }$Q$ \textit{is the largest subgroup of
admissible translators for }$\mu $-\textit{right quasi-invariance:}%
\[
\mu _{g}\sim \mu \text{ for }g\in Q. 
\]%
\textit{In particular, if }$Q_{R}$ \textit{is a subgroup, then }$Q=Q_{R}.$

\bigskip

\noindent \textbf{Proof. }For $K\in \mathcal{K}(G):$ if $\mu (K)=0$ and $%
g\in Q$, then $\mu (Kg)=0.$ For suppose otherwise. Then $\mu (Kg)>0,$ and
then also $\mu ((Kg)g^{-1})>0,$ as $g^{-1}\in Q,$ i.e. $\mu (K)>0,$ a
contradiction. Evidently, if $g$ admits right quasi-invariance, then $\mu
(Kg)>0$ for $K\in \mathcal{K}_{+}(\mu );$ so if $g$ lies in a subgroup
admitting right quasi-invariance, then also $\mu (Kg^{-1})>0,$ and so $g\in
Q.$ $\square $

\bigskip

\noindent \textbf{Proposition 4.4 (Subgroups).} \textit{For }$G$\textit{\
abelian and symmetric }$\mu \in \mathcal{P}(G),$ $Q_{R}$ \textit{is a
subgroup of} $G.$ \textit{In particular, for }$G$ \textit{a Hilbert space as
in \S 1, equipped with a symmetric Gaussian measure }$\mu =\gamma $\textit{,
the Cameron-Martin space }$H(\gamma )$ \textit{is precisely of the form} $%
Q_{R}.$

\bigskip

\noindent \textbf{Proof. }For $\mu $ symmetric, and $K\in \mathcal{K}(G):$
if $\mu (K)>0,$ then $\mu (K^{-1})>0;$ so, for $x\in Q_{R}$ , $\mu
(K^{-1}x)>0.$ But, as $G$ is abelian,%
\[
\mu (Kx^{-1})=\mu (xK^{-1})=\mu (K^{-1}x)>0. 
\]%
So if $x\in Q_{R},$ then $x^{-1}\in Q_{R},$ i.e. $Q_{R}=Q$, and, by Lemma
4.1, $Q_{R}$ is a group.

In particular, for any (symmetric) Gaussian $\mu =\gamma $ with domain a
Hilbert space $G,$ as in \S 1, the Cameron-Martin space $H(\gamma )$
coincides with $Q_{R}(\gamma )$. Indeed, the Hilbert space $G,$ regarded as
an additive group, is abelian, and $h\in H(\gamma )$ holds iff $\mu _{h}$ is
equivalent to $\mu $. This may be re-stated as follows: $h\in H(\gamma )$
holds iff for all $K\in \mathcal{K}(G):$ $\mu (K)>0$ iff $\gamma _{h}(K)>0.$ 
$\square $

\bigskip

\noindent \textbf{Lemma 4.2. }$Q_{R}\mathcal{N}(\mu )\subseteq \mathcal{N}%
(\mu )$ and $Q\mathcal{N}(\mu )=\mathcal{N}(\mu ).$

\bigskip

\noindent \textbf{Proof. }For $h\in Q,$ observe that $\mathcal{K}_{+}(\mu )h=%
\mathcal{K}_{+}(\mu )$: indeed, if $h\in Q,$ then $\mu (Kh)>0$ for all $K\in 
\mathcal{K}_{+}(\mu ),$ so that $\mathcal{K}_{+}(\mu )h\subseteq \mathcal{K}%
_{+}(\mu ).$ Further, for any compact $K\in \mathcal{K}_{+}(\mu ),$ as $%
h^{-1}\in Q$ so that $\mu (Kh^{-1})>0,$ $K=Kh^{-1}h\in \mathcal{K}_{+}(\mu
)h $.

It now follows that if $\mu (Kg)=0$ for all $K\in \mathcal{K}_{+}(\mu ),$
then $\mu (Khg)=0$ for $h\in Q$ and all $K\in \mathcal{K}_{+}(\mu ),$ i.e. $%
hg\in \mathcal{N}(\mu ),$ i.e. $Q\mathcal{N}(\mu )\subseteq \mathcal{N}(\mu
).$ But if $g\in \mathcal{N}(\mu )$ and $h\in Q,$ then likewise $h^{-1}g\in 
\mathcal{N}(\mu )$, as $h^{-1}\in Q,$ so $g=h(h^{-1}g)\in Q\mathcal{N}(\mu
). $ $\square $

\bigskip

\noindent \textbf{Proposition 4.5. }$\mathcal{N}(\mu )$\textit{\ is Borel.}

\bigskip

We give two proofs below, the second of which was kindly suggested by the
Referee.

\bigskip

\noindent \textbf{First Proof. }The sets $\{g:\mu (Kg)<1/n\}$ are open as $%
x\mapsto \mu (Kx)$ is upper semicontinuous for $K$ compact. Indeed, the set%
\[
\{(g,K):\mu (Kg)<1/n\} 
\]%
is open in the product space $G\times \mathcal{K}(G)$: for open $U\supseteq
K\ $with $\mu (U)<1/2n,$ choose $\delta >0$ with $KgB_{2\delta }\subseteq U.$
Then for compact $H\subseteq KB_{\delta }$ and $h\in gB_{\delta },$%
\[
Hh\subseteq KgB_{2\delta }\subseteq U. 
\]

By inner regularity we may assume that $\mu $ is concentrated on a $\sigma $%
-compact set, say on $\dbigcup\nolimits_{n}K_{n}$ with $K_{n}$ an ascending
sequence of compact sets. Then $\{K\subseteq K_{m}:\mu
(K)=0\}=\dbigcap\nolimits_{n\in \mathbb{N}}\{K\subseteq K_{m}:\mu (K)<1/n\}$
is $\mathcal{G}_{\delta }$ in $\mathcal{K}(K_{m}),$ and its complement an $%
\mathcal{F}_{\sigma }$ in $\mathcal{K}(G),$ as $\mathcal{K}(K_{m})$ is
compact. Consequently,%
\[
\{(g,K):K\in \mathcal{K}_{+}(\mu ),\mu (Kg)\geq
1/n\}=\dbigcup\nolimits_{m}\{(g,K):K_{m}\cap K\in \mathcal{K}_{+}(\mu ),\mu
(Kg)\geq 1/n\} 
\]%
\[
=\{(g,K):\mu (Kg)\geq 1/n\}\cap \dbigcup\nolimits_{m}\{(g,K):K_{m}\cap K\in 
\mathcal{K}_{+}(\mu )\}\in \mathcal{F}_{\sigma }(G\times \mathcal{K}(G)). 
\]%
So%
\[
G\backslash \mathcal{N}(\mu )=\mathrm{proj}_{G}\dbigcup\nolimits_{m,n}%
\{(g,K):K_{m}\cap K\in \mathcal{K}_{+}(\mu ),\mu (Kg)\geq 1/n\}. 
\]%
But the vertical sections $\{g\}\times \dbigcup\nolimits_{m}\{K:K_{m}\cap
K\in \mathcal{K}_{+}(\mu ),\mu (Kg)\geq 1/n\}$ are $\sigma $-compact. So the
projection is Borel, by the Arsenin-Kunugui theorem ([Rog, Th. 1.4.3, Th.
5.9.1], cf. [Kec, Th. 18.18]). $\square $

\bigskip

\noindent \textbf{Second Proof. }Assume w.l.o.g. that $\mu \in \mathcal{P}(G)
$ so that also $\mu _{g}\in \mathcal{P}(G)$ for $g\in G.$ (Otherwise, since $%
\mu $ is $\sigma $-finite (\S 2), we may replace $\mu $ by an equivalent
probability measure $\mu ^{\prime }$; then $\mu _{g}\bot \mu $ iff $\mu
_{g}^{\prime }\bot \mu ^{\prime }$.) Since $G$ is separable, we may choose a
(dense) sequence of continuous real-valued functions $\{f_{n}\}$ with $%
||f_{n}||_{\infty }\leq 1$ so that, with $||.||_{TV}$ denoting the total
variation distance, 
\[
||\mu -\mu _{g}||_{TV}=\sup_{n}\int f_{n}d(\mu -\mu _{g})\quad (g\in G);
\]%
cf. [Con, App. C Cor. C.14], [Yos, I.3 Cor.]. Take $\{K_{m}\}_{m}$ an
increasing sequence of compact sets with $\mu $ concentrated on on its
union; then $f_{nm},$ where 
\[
f_{nm}(g):=(f_{n}\ast 1_{K_{m}^{-1}})(g^{-1}),
\]%
is continuous, cf. [HewR, 20.16] or [Rud1, \S 1.1.5-1.1.6]. But%
\[
f_{nm}(g)=\int_{K_{m}}f_{n}(yg^{-1})d\mu (y)\rightarrow \int f_{n}(x)d\mu
_{g}(x)\quad \quad \text{(as }m\rightarrow \infty ).
\]%
So%
\[
||\mu -\mu
_{g}||_{TV}=\sup_{n}\lim_{m}\int_{K_{m}}[f_{n}(y)-f_{n}(yg^{-1})]d\mu (y),
\]%
and so $f(g):=||\mu -\mu _{g}||_{TV}$ ($g\in G)$ is Borel measurable. But,
cf. [Bog1, Problem 2.21], $\mathcal{N}(\mu )=f^{-1}\{2\},$ which is thus
Borel. $\square $

\section{Groups versus vector spaces}

Here we link our new results, in a group context, to classical
Cameron-Martin theory, in a topological vector space context. There is some
similarity here to material in [Hey1, \S 3.4] on (homomorphic) embeddings of 
$\mathbb{Q}$ (rational embeddability) and of $\mathbb{R}$ (the more
exacting, continuous embeddability) in the space of probability measures;
for later developments see [Hey2], [HeyP]. The latter are related to
divisibility properties of groups (and of the convolution semigroups of
measures). Recall that a group $G$ is (infinitely, or $\mathbb{N}$-) \textit{%
divisible} [HewR, A5] if for each $n\in \mathbb{N}$ every element $g\in G$
has an $n$-th root $h\in G$, i.e. with $h^{n}=g$ (for their structure theory
in the abelian case, see e.g. [HewR, A14], [Fuc], or [Kap]).

The refinement norm used to define the classical Cameron-Martin space within
a topological vector space $X$ depends ultimately on the embedding of the
continuous linear functionals, $X^{\ast },$ in $L^{2}(X,\gamma )$ and on the
Riesz Representation theorem for Hilbert spaces. In common with the
additive-subgroup literature of topological vector spaces (for the
corresponding Pontryagin-van Kampen duality theory, see e.g. Gel$^{\prime }$%
fand-Vilenkin [GelV], [Bana], [Mor], and [Xia]), we consider the natural
group analogue here to be the embedding of continuous real-valued additive
maps on a metrizable abelian group $X$ into $L^{2}(X,\mu )$ for $\mu \in 
\mathcal{P}(X),$ with mean (`averaging map')%
\[
a_{\mu }:x^{\ast }\mapsto \mu (x^{\ast }):=\int_{X}x^{\ast }(x)d\mu
(x)\qquad (x^{\ast }\in X^{\ast }), 
\]%
with $X^{\ast }$ the continuous additive maps on $X.$ One may then, as in
the classical setting ([Bog1, 2.2]), define a covariance operator by%
\[
R(x^{\ast })(y^{\ast }):=\int_{X}[x^{\ast }(x)-\mu (x^{\ast })][y^{\ast
}(x)-\mu (y^{\ast })]d\mu (x)\quad \quad (x^{\ast },y^{\ast }\in X^{\ast }). 
\]%
As there, for $h\in X$ define the \textit{Cameron-Martin group-norm}%
\[
|h|_{H}:=\sup \{x^{\ast }(h):x^{\ast }\in X^{\ast },R(x^{\ast },x^{\ast
})\leq 1\}, 
\]%
\[
H:=\{h\in X:|h|_{H}<\infty \}. 
\]%
It may now readily be checked that $H$ is a subgroup, and that $|h|_{H}$ is
a group-norm on $H$. (For any $h\in X,$ with $h\neq 1_{X},$ use a standard
extension theorem, e.g. as in [HewR, A.7], to extend the partial
homomorphism $h^{s}\mapsto s$ (for $s\in \mathbb{Z}$) to a full homomorphism 
$x_{h}^{\ast },$ say of unit variance; then $|h|_{H}\geq x_{h}^{\ast }(h)>0,$
as $x_{h}^{\ast }$ is non-constant. This is an analogue of the Gelfand-Ra%
\u{\i}kov theorem on point separation by characters [HewR, 22.12], cf.
[Tar]. See [Bad, Lemma 1] for a variant Hahn-Banach extension theorem in the
present commutative group context.) It is not clear, however, whether the
resulting group is trivial for $\mu \in \mathcal{P}(X)$. In the classical
Gaussian $\gamma $ context, $|h|_{H}=\infty $ implies the mutual singularity 
$\gamma _{h}\bot \gamma $ [Bog1, 2.4.5(i)].

For $h\in H$ and $n\in \mathbb{N},$ assuming w.l.o.g. that the supremum for $%
|h|_{H}$ occurs with $x^{\ast }(h)>0,$ 
\begin{eqnarray*}
|nh|_{H} &=&\sup \{x^{\ast }(nh):x^{\ast }\in X^{\ast },R(x^{\ast },x^{\ast
})\leq 1\} \\
&=&\sup \{nx^{\ast }(h):x^{\ast }\in X^{\ast },R(x^{\ast },x^{\ast })\leq 1\}
\\
&=&n|h|_{H}.
\end{eqnarray*}%
Thus the norm $|.|_{H}$ (which is subadditive) may be said to be $\mathbb{N}$%
\textit{-homogeneous}, as in [BinO2, \S 3.2], or \textit{sublinear} in the
sense of Berz (see below).

Suppose now that the group $X$ is (infinitely) divisible, so $x/n$ is
defined for $n\in \mathbb{N}$, as is also $qx$ for rational $q$. It follows
by a similar argument to that above that if $h\in H$ and $q>0,$ then $%
|qh|_{H}<\infty ,$ and so $H$\ is also divisible. We now briefly study
sublinear group-norms on a divisible abelian group, proving in particular in
Prop. 5.1 that $H$ is a topological vector space.

\bigskip

\textbf{Definition.} In an abelian ($\mathbb{N}$-) divisible group $G,$
since its group-norm $||.||$ is subadditive, we follow Berz [Ber] (cf.
[BinO3,5,6]) in calling it \textit{sublinear} if%
\[
||ng||=n||g||\qquad (g\in G,n\in \mathbb{N}) 
\]%
(so that $||g/n||=||g||/n),$ or equivalently and more usefully:%
\[
||qg||=q||g||\qquad (g\in G,q\in \mathbb{Q}_{+}). 
\]

\noindent \textbf{Remark.} The triangle inequality for the norm gives $%
||g||\leq n||g/n||;$ the definition requires the reverse inequality.

\bigskip

\noindent \textbf{Proposition 5.1.} \textit{A divisible abelian topological
group complete under a sublinear group-norm is a topological vector space
under the action} 
\begin{equation}
t\cdot g:=\lim_{q\rightarrow t}qg\qquad (t\in \mathbb{R}_{+},\text{ }q\in 
\mathbb{Q}_{+},\text{ }g,\in G).  \tag{$\dag $}
\end{equation}%
\textit{In particular, the group has the embeddability property defined by} $%
t\mapsto t\cdot g,$ \textit{for }$t\in \mathbb{R}_{+}$.

\bigskip

\noindent \textbf{Proof. }Since 
\[
||(q-q^{\prime })(g-g^{\prime })||=|q-q^{\prime }|\cdot ||g-g^{\prime
}||\qquad (q,q^{\prime }\in \mathbb{Q}_{+},\text{ }g,g^{\prime }\in G), 
\]%
and the group is norm-complete, the action ($\dag $) is well defined. It is
also jointly continuous, since passing to appropriate limits over $\mathbb{Q}%
_{+}$%
\[
||(t-t^{\prime })(g-g^{\prime })||=|t-t^{\prime }|\cdot ||g-g^{\prime
}||\qquad (t,t^{\prime }\in \mathbb{R}_{+},\text{ }g,g^{\prime }\in G). 
\]%
The action extends to $\mathbb{R}$ by taking $-t\cdot g:=t\cdot (-g),$ and
converts the group into a vector space, as $(q+q^{\prime })g=qg+q^{\prime }g$
and $q(g+g^{\prime })=qg+qg^{\prime }$ may be taken to the limit through $%
\mathbb{Q}$. Likewise the final statement follows, as%
\[
||t\cdot g||=\lim_{q\rightarrow t}||qg||=t||g||.\qquad \square 
\]

\bigskip

We verify that the classical boundedness theorem holds in the group context
for additive functions, and hence that such functions are linear (in the
sense of the action ($\dag $)).

\bigskip

\noindent \textbf{Proposition 5.2.} \textit{For an abelian divisible group }$%
G$\textit{\ complete under a sublinear group-norm and} $x^{\ast
}:G\rightarrow \mathbb{R}$ \textit{additive, }$x^{\ast }$\textit{\ is
continuous iff}%
\[
||x^{\ast }||=\sup \{|x^{\ast }(x)|/||x||:x\in G\backslash \{0\}\}<\infty , 
\]%
\textit{and then}%
\[
|x^{\ast }(x)|\leq ||x^{\ast }||\cdot ||x||. 
\]

\bigskip

\noindent \textbf{Proof.} If $x^{\ast }$ is continuous, choose $\delta >0$
with $|x^{\ast }(x)|\leq 1$ for $||x||\leq \delta .$ For $x\neq 0$ and
rational $q\geq ||x||/\delta ,$ as $||x/q||=||x||/q\leq \delta ,$%
\[
|x^{\ast }(x)|=q|x^{\ast }(x/q)|\leq q; 
\]%
then taking limits as $q\rightarrow ||x||/\delta $ through $\mathbb{Q}$
yields%
\[
|x^{\ast }(x)|\leq ||x||/\delta . 
\]%
This holds also for $x=0.$ So $||x^{\ast }||\leq 1/\delta <\infty .$ In this
case, by definition, $|x^{\ast }(x)|\leq ||x^{\ast }||\cdot ||x||$ for $%
x\neq 0,$ and this again holds also for $x=0.$

Conversely, if $||x^{\ast }||<\infty ,$ then again $|x^{\ast }(x)|\leq
||x^{\ast }||\cdot ||x||$ , and so $x^{\ast }$ is continuous at $0$ and
hence everywhere. $\square $

\bigskip

\noindent \textbf{Corollary 5.1.}\textit{\ A continuous additive function }$%
x^{\ast }$\textit{\ may be extended by taking}%
\[
x^{\ast }(t\cdot g):=\lim_{q\rightarrow t}x^{\ast }(qg)\qquad (g\in G,\text{ 
}q\in \mathbb{Q}_{+},\text{ }t\in \mathbb{R}_{+}), 
\]%
\textit{which then makes it linear in the sense of the action} $t\cdot g:$%
\[
x^{\ast }(t\cdot g)=t\cdot x^{\ast }(g)\qquad (g\in G,t\in \mathbb{R}_{+}). 
\]%
\noindent \textbf{Proof. }Convergence of\textbf{\ }$x^{\ast }(qg)$ as $%
q\rightarrow t$ follows from the finiteness of $||x^{\ast }||$ and 
\[
|x^{\ast }(qg)-x^{\ast }(q^{\prime }g)|=|x^{\ast }((q-q^{\prime })g)|\leq
||x^{\ast }||\cdot |q-q^{\prime }|\cdot ||g||. 
\]%
The final claim follows, taking limits through $\mathbb{Q}_{+},$%
\[
x^{\ast }(tg)=\lim_{q\rightarrow t}x^{\ast }(qg)=\lim_{q\rightarrow t}q\cdot
x^{\ast }(g)=t\cdot x^{\ast }(g).\qquad \square 
\]

\bigskip

Returning to the context of the Cameron-Martin group-norm, notwithstanding
the issue of possible triviality of $H$, one may define here as in the
classical context the set $X_{\mu }^{\ast }\subseteq L^{2}(\mu ),$ by
passing to the closed span in $L^{2}(\mu )$ of the centred (mean-zero) image 
$\{x^{\ast }-\mu (x^{\ast }):x^{\ast }\in X^{\ast }\}$ of $X^{\ast }$. Then
one may extend the domain of $R$ above to the subspace $X_{\mu }^{\ast }$ of 
$L^{2}(\mu ),$ so that $R(f)$ for $f\in X_{\mu }^{\ast }$ is given by 
\[
R(f)(y^{\ast }):=\int f(x)[y^{\ast }(x)-\mu (y^{\ast })]d\mu (x)\qquad
(y^{\ast }\in X^{\ast }). 
\]%
(For $f=x^{\ast }\in X^{\ast },$ this coincides with the previous
definition, since $\mu (x^{\ast })\int_{X}[y^{\ast }(x)-\mu (y^{\ast })]d\mu
(x)=0.$) Then $h\in H$ if $h=R(g)$ for some $g=\hat{h}\in X_{\mu }^{\ast }$;
then $\delta _{h}(y^{\ast })=y^{\ast }(h)=R(g)(y^{\ast }),$ for $y^{\ast
}=z^{\ast }-\mu (z^{\ast })$ with $z^{\ast }\in X^{\ast }$, and so taking $%
y^{\ast }\rightarrow g$%
\[
|h|_{H}=||g||_{L^{2}(\mu )}. 
\]%
Conversely, for fixed $h\in H,$ the map $f\mapsto f(h)$, for $f\in X_{\mu
}^{\ast }$ (including $f=x^{\ast }-\mu (x^{\ast })$ for $x^{\ast }\in
X^{\ast }$) is represented by%
\[
\langle f,\hat{h}\rangle _{L^{2}(\mu )}, 
\]%
by the Riesz Representation Theorem. So $h=R(\hat{h}).$

For $h,k\in H,$ define their inner product by referring to the elements $%
\hat{h},\hat{k}$ such that $h=R(\hat{h}),$ $k=R(\hat{k}),$ passing in $R\ $%
to the limit, and putting%
\[
(h,k)_{H}:=\int_{X}\hat{h}(x)\hat{k}(x)d\mu (x). 
\]%
Then $H\ $inherits an inner-product structure by%
\[
|h|_{H}^{2}=||\hat{h}||_{L^{2}(\mu )}^{2}=\int_{X}(\hat{h}(x))^{2}d\mu
(x)=(h,h)_{H(\mu )}. 
\]%
In the vector-space Gaussian context, for $h\in H$ the density of $\mu _{h}$
w.r.t. $\mu $ is given by $(CM)$ above in \S 1.

For Gaussian measures on (locally compact) groups $G$, see \S 9.12.

\section{Reference measures with selective subcontinuity}

We briefly recall the construction of a measure due to Solecki, as this
exemplifies a property at the heart of the construction of a Steinhaus
triple in \S 7. For this it will be helpful to review from Th. 3.4 the
notation $\mu _{\_}^{\mathbf{t}}$ for $\mathbf{t}=\{t_{n}\}$ a null
sequence, i.e. with $t_{n}\rightarrow 1_{G}:$ 
\[
\mu _{-}^{\mathbf{t}}(K):=\lim \inf\nolimits_{n\rightarrow \infty }\mu
(Kt_{n})\qquad (K\in \mathcal{K}(G)). 
\]%
When $\mu (K)=\mu _{-}^{\mathbf{t}}(K)$ we refer to \textit{selective
subcontinuity} of $\mu $ `along' $\mathbf{t}$ (cf. [BinO7]). Thus, for $G$
not locally compact and $K\in \mathcal{K}_{+}(\mu )$, by Cor. 3.2 there will
be $\mathbf{t}$ with $\mu (Kt_{n})\rightarrow 0,$ but there may, and in
certain specified circumstances necessarily will, also exist $\mathbf{t}$
with $\mu _{-}^{\mathbf{t}}(K)>0$ (as when the Subcontinuity Theorem, Th.
6.1 below, holds).

Recall that a (Polish) group $G$ is \textit{amenable at} $1$ ([Sol]; cf.
[BinO7, \S 2] for the origin of this term) if, given any sequence of
measures $\nu _{n}\in \mathcal{P}(G)$ with $1_{G}\in $ \textrm{supp}$(\nu
_{n}),$ there are $\sigma $ and $\sigma _{n}$ in $\mathcal{P}(G)$ with $%
\sigma _{n}\ll \nu _{n}$ for each $n\in \mathbb{N}$, and 
\[
\sigma _{n}\ast \sigma (K)\rightarrow \sigma (K)\qquad (K\in \mathcal{K}%
(G)). 
\]%
Here $\ast $ denotes convolution. (Abelian Polish groups all have this
property [Sol, Th. 3.2].) With $\delta _{g}$ the Dirac measure at $g$ and $%
\mathbf{t}$ a non-trivial null sequence, taking 
\[
\nu _{n}:=2^{n-1}\sum\nolimits_{m\geq n}2^{-m}\delta _{t_{m}^{-1}}\in 
\mathcal{P}(G) 
\]%
$,$ we denote by $\sigma _{n}(\mathbf{t})$ and $\sigma (\mathbf{t})$ the
measures whose existence the definition above asserts. We term $\sigma (%
\mathbf{t})$ a \textit{Solecki reference measure}. Below we place further
restrictions on the rate of convergence of $\mathbf{t}$ and refer to a
symmetrized version of $\nu _{n}.$

In the above setting, $\sigma (K)>0$ implies that for some sequence $%
m(n):=m_{K}(n)$%
\[
\sigma (Kt_{m_{K}(n)})\rightarrow \sigma (K). 
\]%
We repeat here a result from the companion paper [BinO7], as its proof
motivates the definition which follows.

\bigskip

\noindent \textbf{Theorem 6.1 }(\textbf{Subcontinuity Theorem}, after
Solecki [Sol, Th. 1(ii)]).\textit{\ For }$G\ $\textit{amenable at }$1$%
\textit{, }$0<\theta <1,$ \textit{and }$\mathbf{t}$ \textit{a null sequence,
there is }$\sigma =\sigma (\mathbf{t})\in \mathcal{P}(G)$ \textit{such that
for each }$g\in G,K\in \mathcal{K}(G)$ \textit{with }$\sigma (gK)>0$\textit{%
\ there is a subsequence }$\mathbf{s}=\mathbf{s}(g,K):=\{t_{m(n)}\}$ \textit{%
with}%
\[
\sigma (gKt_{m(n)})>\theta \sigma (gK)\qquad \text{(}n\in \mathbb{N}\text{)},%
\text{ so}\qquad \sigma _{\_}^{\mathbf{s}}(gK)>0. 
\]%
\textit{That is, }$\sigma $\textit{\ is subcontinuous along }$\mathbf{s}$%
\textit{\ on }$gK.$ \textit{In particular}%
\[
\sigma (K)=\sigma _{\_}^{\mathbf{t}}(K). 
\]

\noindent \textbf{Proof.} For $\mathbf{t}=\{t_{n}\}$ null, put $\nu
_{n}:=2^{n-1}\sum\nolimits_{m\geq n}2^{-m}\delta _{t_{m}^{-1}}\in \mathcal{P}%
(G);$ then $1_{G}$ $\in $ \textrm{supp}$(\nu _{n})\supseteq
\{t_{m}^{-1}:m>n\}.$ By definition of amenability at $1$, in $\mathcal{P}(G)$
there are $\sigma $ and $\sigma _{n}\ll \nu _{n},$ with $\sigma _{n}\ast
\sigma (K)\rightarrow \sigma (K)$ for all $K\in \mathcal{K}(G).$

Fix $K\in \mathcal{K}(G)$ and $g$ with $\sigma (gK)>0.$ As $gK$ is compact, $%
\sigma _{n}\ast \sigma (gK)\rightarrow \sigma (gK);$ then w.l.o.g. 
\begin{equation}
\nu _{n}\ast \sigma (gK)>\theta \sigma (gK)\qquad (n\in \mathbb{N}). 
\tag{$\ddagger $}
\end{equation}%
For $m,n\in \mathbb{N}$ choose $\alpha _{mn}\geq 0$ with $%
\sum\nolimits_{m\geq n}\alpha _{mn}=1$ $(n\in \mathbb{N})$ and $\nu
_{n}:=\sum\nolimits_{m\geq n}\alpha _{mn}\delta _{t_{m}^{-1}}.$ Then for
each $n$ there is $m=m(n)\geq n$ with%
\[
\sigma (gKt_{m})>\theta \sigma (gK); 
\]%
otherwise, summing the reverse inequalities over all $m\geq n$ contradicts ($%
\ddagger $). So $\lim_{n}\sigma (gKt_{m(n)})\geq \theta \sigma (gK):$ $%
\sigma $ is subcontinuous along $\mathbf{s}$ $:=\{t_{m(n)}\}$ on $gK.$

For the last assertion, take $g=1_{G}$ and recall that $\sigma (K)\geq
\sigma _{\_}^{\mathbf{t}}(K),$ by upper semi-continuity of $t\mapsto \sigma
(Kt).$ $\square $

\bigskip

A corollary to this is a `non-Haar-null' version of the Steinhaus-Weil
interior point theorem:

\bigskip

\noindent \textbf{Theorem SW} (\textbf{Steinhaus-Weil Theorem beyond Haar }%
[Sol, Th 1(ii)]).\textit{\ For }$G\ $\textit{amenable at }$1:$\textit{\ if} $%
E\in \mathcal{U}(G)$\textit{\ is not left Haar null, then }$1_{G}\in \mathrm{%
int}(E^{-1}E).$

\bigskip

\noindent \textbf{Proof.} Suppose otherwise; then $1_{G}\notin \mathrm{int}%
(E^{-1}E)$ and we may choose $t_{n}\in B_{1/n}\backslash E^{-1}E.$ As $%
t_{n}\rightarrow 1_{G},$ choose $\sigma =\sigma (\mathbf{t})$ as in the
preceeding theorem. Since $E$ is not left Haar null, $\sigma (gE)>0$ for
some $g$. For this $g,$ choose compact $K\subseteq E$ with $\sigma (gK)>0.$
Then by the Subcontinuity Theorem (Th. 6.1) and by Prop. 3.1 for $\Delta
=\sigma _{-}^{\mathbf{t}}(K)/4,$ there is $\delta >0$ such that%
\[
\emptyset \neq B_{\delta }^{\Delta }\subseteq K^{-1}g^{-1}gK=K^{-1}K; 
\]%
moreover, as in Prop. 3.1, $t_{n}\in K^{-1}K\subseteq E^{-1}E$ for
infinitely many $n$, which contradicts the choice of $\mathbf{t}.$ So $%
1_{G}\in \mathrm{int}(E^{-1}E).$ $\square $

\bigskip

\noindent \textbf{Definition} ([BinO7]). Say that a null sequence $\mathbf{t}
$ is \textit{regular} if $\mathbf{t}$ is non-trivial, $\{||t_{k}||\}_{k}$ is
non-increasing, and 
\[
||t_{k}||\leq r(k):=1\left/ [2^{k}(k+1)]\right. \qquad (k\in \mathbb{N}). 
\]%
For regular $\mathbf{t}$, put 
\[
\tilde{\nu}_{k}=\tilde{\nu}_{k}(\mathbf{t}):=2^{k-1}\sum\nolimits_{m\geq
k}2^{-m}(\delta _{t_{m}^{-1}}+\delta _{t_{m}})=\frac{1}{2}\delta
_{t_{k}^{-1}}+\frac{1}{4}\delta _{t_{k+1}^{-1}}+.... 
\]%
Then $\tilde{\nu}_{k}(B_{r(k)})=1$ and $\tilde{\nu}_{k}$ is \textit{symmetric%
}: $\tilde{\nu}_{k}(K^{-1})=\tilde{\nu}_{k}(K),$ as $t\in K$ iff $t^{-1}\in
K^{-1}.$ In the definition below, which is motivated by the proof of Theorem
6.1, we will view the measure $\tilde{\nu}_{k}$ as just another version of $%
\nu _{k}$ above (by merging $\mathbf{t}^{-1}$ with $\mathbf{t}$ by
alternation of terms).

\bigskip

\noindent \textbf{Definition} ([BinO7]). Say that a (Polish) group $G$ is 
\textit{strongly amenable at 1} if $G$ is amenable at 1, and for each 
\textit{regular} $\mathbf{t}$ the Solecki measure $\sigma (\mathbf{t})$
corresponding to $\nu _{k}(\mathbf{t})$ has associated measures $\sigma _{k}(%
\mathbf{t})\ll \nu _{k}(\mathbf{t})$ with the following \textit{%
concentration property}. Writing, for $k\in \mathbb{N}$, 
\[
\sigma _{k}:=\sum\nolimits_{m\geq k}a_{km}\delta _{t_{m}^{-1}}, 
\]%
for some non-negative sequences $\mathbf{a}_{k}:=%
\{a_{kk,}a_{k,k+1,}a_{k,k+2},...\}$ of unit $\ell _{1}$-norm, there is an
index $j$ and $\alpha >0$ with 
\[
a_{k,k+j}\geq \alpha >0\qquad \text{for all large }k. 
\]

A refinement of Solecki's proof of the Subcontinuity Theorem (Th. 6.1 above)
yields the following two results, for the proof of which we refer the reader
to the companion paper [BinO7].

\bigskip

\noindent \textbf{Theorem 6.2 }(\textbf{Strong amenability at 1}, [BinO7,
Th. 4] after [Sol2, Prop. 3.3(i)]). \textit{Any abelian Polish group }$G$%
\textit{\ is strongly amenable at 1.}

\bigskip

\noindent \textbf{Theorem 6.1}$_{\text{\textbf{S}}}$ \textbf{(Strong
Subcontinuity Theorem). }\textit{For }$G$\textit{\ a (Polish) group that is
strongly amenable at 1, if }$\mathbf{t}$\textit{\ is regular and }$\sigma
=\sigma (\mathbf{t})$\textit{\ is a Solecki measure -- then for }$K\in 
\mathcal{K}_{+}(\sigma )$%
\[
\sigma (K)=\lim_{n}\sigma (Kt_{n})=\sigma _{-}^{\mathbf{t}}(K). 
\]

\noindent \textbf{Remark. }The reference measure $\sigma (\mathbf{t})$ in
the last theorem may in fact be selected symmetric [BinO7], in which case%
\textit{\ }$Q_{R}(\sigma (\mathbf{t}))$ is a subgroup.

\bigskip

We note an immediate corollary, useful in \S 8 below.

\bigskip

\noindent \textbf{Corollary 6.1}. \textit{For }$G,\mathbf{t}$\textit{\ and }$%
\sigma $ \textit{as in Th. 6.1}$_{\text{S}}$ \textit{above, and }$K,H\in 
\mathcal{K}(G),\delta >0$\textit{: if }$0<\Delta <\sigma (K)$\textit{\ and }$%
0<D<\sigma (H),$\textit{\ then there is }$n$\textit{\ with} 
\[
B_{\delta }^{K\Delta }\cap B_{\delta }^{HD}\supseteq \{t_{m}:m\geq n\}. 
\]

\bigskip

\noindent \textbf{Proof.} Take $\varepsilon :=\min \{\sigma (K)-\Delta
,\sigma (H)-D\}>0.$ As $K,H\in \mathcal{K}_{+}(\sigma ),$ there is $n$ such
that $||t_{m}||<\delta $ for $m\geq n$ and 
\[
\sigma (Kt_{m})\geq \sigma (K)-\varepsilon \geq \Delta ,\qquad \sigma
(Ht_{m})\geq \sigma (H)-\varepsilon \geq D\qquad (m\geq n).\qquad \square 
\]

\section{The Steinhaus support{\protect\large \ }of a measure}

In this section we construct (one might say via a `disaggregation') the
Steinhaus support $H(\mu )$ of a probability measure $\mu $ defined on a
Polish group $G$ (see Th. 7.1); this is possible provided the measure has
`sufficient subcontinuity' (defined below) -- sufficient to allow a
relative-interior Steinhaus property, relative to some embedded `subspace'.
In the next section we apply the construction to a Solecki measure $\sigma (%
\mathbf{t})$ for a regular null sequence $\mathbf{t}$ as in \S 6.

\bigskip

\textbf{Definition.} Say that a probability measure has \textit{sufficient
subcontinuity}, written $\mu \in \mathcal{P}_{\text{suf}}(G),$ if for all $%
K\in \mathcal{K}_{+}(\mu )$ and $\delta >0$ there is $\Delta (K,\delta )\geq
0$ so small that for $\Delta (K,\delta )<\Delta <\mu (K)$%
\[
B_{\delta }^{K,\Delta }=B_{\delta }^{K,\Delta }(\mu )=\{s\in B_{\delta }:\mu
(Ks)>\Delta \} 
\]%
is infinite. Above, if $\Delta (K,\delta )\equiv 0,$ say that $\mu $ is of 
\textit{Solecki type}; these will be considered in \S 8.

\bigskip

Lemma 7.1 below asserts that, for $G$ amenable at $1,$ a Solecki measure $%
\mu =\sigma (\mathbf{t})$ has this property. Further motivation for working
with this assumption is provided by Th 3.6: for $(H,G,\mu )$ a Steinhaus
triple and $K\in \mathcal{K}_{+}(\mu ),$ the set $\mathcal{O}(K):=\{s\in
H:\mu (Ks)>0\}$ is open in $H$ and $\{s\in H:\mu (Ks)=\mu _{\_}^{H}(Ks)>0\}$
is dense in $\mathcal{O}(K)$.

The goal here is to create a new topology, if not on $G$ then on a dense
subset of $G,$ in which the sets 
\[
B_{\delta }^{K,\Delta }:=\{z\in B_{\delta }:\mu (Kz)>\Delta \} 
\]%
(with $\mu $ understood from context) shall be open sub-base members for
selected $K\in \mathcal{K}_{+}(\mu )$. This is tantamount to requiring that
the corresponding maps $z\mapsto \mu (Kz)$ be continuous on some subset of $%
G;$ cf. Th. 3.2$^{\prime }$ , also by way of justification.

\bigskip

Whenever we consider sets $B_{\delta }^{K,\Delta }$ for $K\in \mathcal{K}%
_{+}(\mu )$ and $\delta >0$ we implicitly assume that $\Delta \geq \Delta
(K,\delta ).$

Notice the monotonicities: 
\[
\Delta \leq D\Longrightarrow B_{\delta }^{K,D}\subseteq B_{\delta
}^{K,\Delta },\qquad \eta \leq \varepsilon \Longrightarrow B_{\eta
}^{K,\Delta }\subseteq B_{\varepsilon }^{K,\Delta }. 
\]

\bigskip

The sequence of Lemmas 7.1-7.4 below justifies the introduction of a new
topology with sub-basic sets of the form $gB_{\delta }^{K,\Delta },$ but 
\textit{only} on those points of $G$ that can be \textit{covered} by these
sets: the detailed statement is in Theorem 7.1 below. The proof strategy
demands both a countable iteration -- an inductive generation of a family of
sets $B_{\delta }^{K,\Delta }$ -- and then a countable subgroup of
translators $g.$ In the subsequent section, we identify which are the points
that can be covered.

For $\mu \in \mathcal{P}(G),\mathcal{H\subseteq K}_{+}(\mu ),$ we put%
\[
\mathcal{B}_{1_{G}}(\mu ,\mathcal{H}):=\{B_{\delta }^{K,\Delta }(\mu ):K\in 
\mathcal{H};\text{ }\delta ,\Delta \in \mathbb{Q}_{+};\mathit{\ }0<\Delta
<\mu (K)\}; 
\]%
this is to be a neighbourhood base at $1_{G}.$ For $\mathcal{H=K}_{+}(\mu )$
we abreviate this to $\mathcal{B}_{1_{G}}(\mu )$ or even to $\mathcal{B}%
_{1_{G}}$, when $\mu $ is understood.

\bigskip

\noindent \textbf{Lemma 7.1. }\textit{For }$\mu \in \mathcal{P}(G)$\textit{, 
}$\mathbf{t}$\textbf{\ }\textit{null and non-trivial, and arbitrary }$\delta
>0$\textit{: if }$0<\Delta <\mu _{-}^{\mathbf{t}}(K),$\textit{\ then }$%
B_{\delta }^{K,\Delta }(\mu )\neq \{1_{G}\}.$\textit{\ In particular, for }$%
G $\textit{\ amenable at }$1$ \textit{and} $\mu =\sigma (\mathbf{t})$\textit{%
: if} $0<\Delta <\sigma (K),$\textit{\ then} $B_{\delta }^{K,\Delta }(\sigma
)$ \textit{is infinite.}

\bigskip

\noindent \textbf{Proof. }Since $\mathbf{t}$ is null and non-trivial, for
all large enough $n$ both\textbf{\ }$t_{n}\in B_{\delta }$ and also $\mu
(Kt_{n})>\Delta .$ For $\mu =\sigma (\mathbf{t})$ and $0<\Delta <\sigma (K),$
pick $0<\theta <1$ with%
\[
\theta \sigma (K)=\Delta . 
\]%
Then for some, necessarily non-trivial, subsequence $\mathbf{s}:=\{s_{n}\}$
of $\mathbf{t}$,%
\[
\sigma (Ks_{n})>\theta \sigma (K)=\Delta . 
\]%
So $B_{\varepsilon }^{K,\Delta }=\{s\in B_{\varepsilon }:\sigma (Ks)>\Delta
\}$ is infinite. $\square $

\bigskip

\noindent \textbf{Lemma 7.2. }\textit{For }$\mu \in \mathcal{P}_{\text{suf}%
}(G)$\textit{\ and }$K\in \mathcal{K}_{+}(\mu ):$\textit{\ if }$w\in
B_{\delta }^{K,\Delta },$\textit{\ then for }$H=Kw$\textit{\ and some }$%
\varepsilon >0$ 
\[
\{w\}\neq wB_{\varepsilon }^{H,\Delta }\subseteq B_{\delta }^{K,\Delta }. 
\]%
\textit{In particular:\newline
}(i)\textit{\ if }$1_{G}\in gB$ \textit{for some }$B\in \mathcal{B}_{1_{G}},$%
\textit{\ then there is is} $B^{\prime }\in \mathcal{B}_{1_{G}}$ \textit{%
with }$1_{G}\in B^{\prime }\subseteq gB;$\newline
(ii) \textit{for }$G$\textit{\ amenable at }$1$\textit{: if} $\mu =\sigma (%
\mathbf{t})$\textit{\ with }$\mathbf{t}$\textit{\ null, then }$%
B_{\varepsilon }^{H,\Delta }$ \textit{may be selected infinite}.

\bigskip

\noindent \textbf{Proof.} As $w\in B_{\delta }$ there is $0<\varepsilon
<\delta $ with $wB_{\varepsilon }\subset B_{\delta }.$ As $w\in B_{\delta
}^{K,\Delta },$ $\mu (H)=\mu (Kw)>\Delta ,$ so, passing to a smaller $%
\varepsilon $ if necessary, there is $\Delta (H,\varepsilon )>0$ so that $%
B_{\varepsilon }^{H,\Delta ^{\prime }}$ is infinite for $\Delta ^{\prime
}\geq \max \{\Delta (H,\varepsilon ),\Delta \}>\Delta (K,\delta )$. Then 
\begin{eqnarray*}
w &\in &wB_{\varepsilon }^{H,\Delta ^{\prime }}=w\{s\in B_{\varepsilon }:\mu
(Hs)>\Delta ^{\prime }\}=\{ws\in wB_{\varepsilon }:\mu (Kws)>\Delta ^{\prime
}\} \\
&\subseteq &\{x\in B_{\delta }:\mu (Kx)>\Delta \}=B_{\delta }^{K,\Delta }.
\end{eqnarray*}

For the last part, suppose $1_{G}\in gB$ with $B=B_{\delta }^{K,\Delta }\in 
\mathcal{B}_{1_{G}};$ then $w\in B$ for $w=g^{-1}.$ Applying the first part,
take $B^{\prime }:=B_{\varepsilon }^{H,\Delta }\in \mathcal{B}_{1}$ for $%
H=Kw $ and the $\varepsilon >0$ above; then, 
\[
w\in wB^{\prime }=B_{\varepsilon }^{H,\Delta }\subseteq B_{\delta
}^{K,\Delta }=B:\qquad 1_{G}\in B^{\prime }\subseteq gB.\qquad \square 
\]

\bigskip

\noindent \textbf{Corollary 7.1.} \textit{If }$x\in yB\cap zC$\textit{\ for }%
$x,y,z\in G$ \textit{and some }$B,C\in \mathcal{B}_{1_{G}},$ \textit{then }$%
x\in x(B^{\prime }\cap C^{\prime })\subseteq yB\cap zC$\textit{\ for some} $%
B^{\prime },C^{\prime }\in \mathcal{B}_{1_{G}}.$

\bigskip

\noindent \textbf{Proof.} As $1_{G}\in x^{-1}yB$ and $1_{G}\in x^{-1}zC,$
there are $B^{\prime },C^{\prime }\in \mathcal{B}_{1_{G}}$ with $1_{G}\in
B^{\prime }\subseteq x^{-1}yB$ and $1_{G}\in C^{\prime }\subseteq x^{-1}zC.$
Then $x\in xB^{\prime }\cap xC^{\prime }\subseteq yB\cap zC.$ $\square $

\bigskip

We now improve on Lemma 7.2 by including some technicalities, whose purpose
is to introduce a \textit{separable} topology on a subspace of $G$ refining
that induced by $\tau _{G}$. In view of the monotonicities observed above,
we may restrict attention to $\delta ,\Delta \in \mathbb{Q}_{+}:=\mathbb{Q}%
\cap (0,\infty ).$

\bigskip

\noindent \textbf{Lemma 7.3. }\textit{For }$\mu \in \mathcal{P}_{\text{suf}%
}(G)$\textit{\ and countable }$\mathcal{H\subseteq K}_{+}(\mu )$\textit{,
there is a countable }$D=D(\mathcal{H})\subseteq G$ \textit{accumulating at }%
$1_{G}$ \textit{such that: if }$w\in B_{\delta }^{K,\Delta }$\textit{\ with }%
$K\in \mathcal{H}$, $\delta ,\Delta \in \mathbb{Q}_{+}$\textit{and }$\Delta
<\mu (K),$ \textit{then for some }$g\in D$\textit{\ with }$\mu (Kg)>\Delta $ 
\textit{and some }$\varepsilon \in \mathbb{Q}_{+},$ 
\[
w\in gB_{\varepsilon }^{Kg,\Delta }\subseteq B_{\delta }^{K,\Delta }. 
\]%
\noindent \textbf{Proof. }As $G$ is separable, we may choose $\{\bar{g}%
_{m}\}_{m\in \mathbb{N}}=\{\bar{g}_{m}(B_{\delta }^{K,\Delta })\}_{m\in 
\mathbb{N}}\subseteq B_{\delta }^{K,\Delta }$ dense in $B_{\delta
}^{K,\Delta }$, an infinite set, by Lemma 7.1. Take%
\[
D=D(\mathcal{H}):=\{\bar{g}_{m}(B_{\delta }^{K,\Delta }):K\in \mathcal{H}%
,\delta ,\Delta \in \mathbb{Q}_{+},\Delta <\mu (K)\}, 
\]%
which is countable. Since $B_{\delta }^{K,\Delta }\subseteq B_{\delta },$ $D$
accumulates at $1_{G}.$ We claim that $D\ $above satisfies the conclusions
of the Lemma.

Fix $w\in B_{\delta }^{K,\Delta },$ with $K,\Delta ,\delta $ as in the
hypotheses. Choose $\varepsilon \in \mathbb{Q}_{+}$ with%
\[
wB_{3\varepsilon }\subseteq B_{\delta }. 
\]%
Choose $\bar{g}_{m}=\bar{g}_{m}(B_{\delta }^{K,\Delta })$ with $||\bar{g}%
_{m}^{-1}w||<\varepsilon ,$ possible by construction of $\{\bar{g}%
_{m}(B_{\delta }^{K,\Delta })\}$. Put $z_{m}:=\bar{g}_{m}^{-1}w;$ then $w=%
\bar{g}_{m}z_{m},$ $z_{m}\in B_{\varepsilon }$ and $\bar{g}_{m}\in
wB_{\varepsilon },$ so $w\in \bar{g}_{m}z_{m}B_{\varepsilon }\in \bar{g}%
_{m}B_{2\varepsilon }\subseteq wB_{3\varepsilon }\subseteq B_{\delta }.$ By
choice of $\{\bar{g}_{m}\}_{m\in \mathbb{N}}$, $\mu (K\bar{g}_{m})>\Delta ,$
and furthermore%
\begin{eqnarray*}
w &\in &\bar{g}_{m}z_{m}\{s\in B_{\varepsilon }:\mu (K\bar{g}%
_{m}z_{m}s)>\Delta \}\subseteq \bar{g}_{m}\{t\in B_{2\varepsilon }:\mu (K%
\bar{g}_{m}t)>\Delta \} \\
&=&\bar{g}_{m}B_{2\varepsilon }^{K\bar{g}_{m},\Delta }\subseteq \{x\in
B_{\delta }:\mu (Kx)>\Delta \}=B_{\delta }^{K,\Delta },
\end{eqnarray*}%
as $\bar{g}_{m}B_{2\varepsilon }\subseteq B_{\delta }.$ $\square $

\bigskip

In Lemma 7.3 above $Kg$ need not belong to $\mathcal{H}$. Lemma 7.4 below
asserts that Lemma 7.3 holds on a countable family $\mathcal{H}$ of compact
sets that is closed under the appropriate translations.

\bigskip

\noindent \textbf{Lemma 7.4.}\textit{\ For }$\mu \in \mathcal{P}_{\text{suf}%
}(G),$ \textit{there are a countable} $\mathcal{H\subseteq K}_{+}(\mu )$ 
\textit{and a countable set} $D=D(\mathcal{H})\subseteq G$ \textit{dense in }%
$G$ \textit{such that: if }$w\in B_{\delta }^{K,\Delta }$\textit{\ with }$%
K\in \mathcal{H}$, $\delta ,\Delta \in \mathbb{Q}_{+}$\textit{and }$0<\Delta
<\mu (K),$ \textit{then for some }$g\in D$\textit{\ with }$\mu (Kg)>\Delta $ 
\textit{with }$Kg\in \mathcal{H}$ \textit{and some }$\varepsilon \in \mathbb{%
Q}_{+},$ 
\[
w\in gB_{\varepsilon }^{Kg,\Delta }\subseteq B_{\delta }^{K,\Delta }. 
\]%
\noindent \textbf{Proof. }Suppose $\mu $ is concentrated on $%
\dbigcup\nolimits_{n}K_{n}$, with the $K_{n}$ compact and $\mu (K_{n})>0$.
Taking $\mathcal{H}_{0}$ to comprise the basic compacts $K_{n}\cap g_{m}\bar{%
B}_{\delta }$ with $\{g_{m}\}$ dense in $G$ and $\delta \in \mathbb{Q}_{+},$
proceed by induction:%
\[
\mathcal{H}_{n+1}:=\{Kg:K\in \mathcal{H}_{n},g\in D(\mathcal{H}_{n}),\delta
,\Delta \in \mathbb{Q}_{+},\mathit{\ }0<\Delta <\mu (K)\}, 
\]%
\[
\mathcal{H}:=\dbigcup\nolimits_{n}\mathcal{H}_{n},\qquad
D:=\dbigcup\nolimits_{n}D(\mathcal{H}_{n}).\qquad \square 
\]

\noindent \textbf{Theorem 7.1. }\textit{For} $\mu \in \mathcal{P}_{\text{suf}%
}(G)$\textit{\ there are a countable} $\mathcal{H\subseteq K}_{+}(\mu )$ 
\textit{and a countable set} $\Gamma =\Gamma (\mathcal{H})\subseteq G$ 
\textit{dense in }$G$ \textit{such that, taking}%
\[
\mathcal{B}_{\mathcal{H}}(\mu )=\{B_{\delta }^{K,\Delta }(\mu )\in \mathcal{B%
}_{1_{G}}:K\in \mathcal{H};\delta ,\Delta \in \mathbb{Q}_{+};\mathit{\ }%
0<\mu (K)<\Delta \}, 
\]%
\[
\mathcal{B}(\mu )=\mathcal{B}_{\Gamma }(\mu ):=\Gamma \cdot \mathcal{B}_{%
\mathcal{H}}(\mu )=\{\gamma B:\gamma \in \Gamma ,B\in \mathcal{B}_{\mathcal{H%
}}(\mu )\} 
\]%
\textit{is a sub-base for a second-countable topology on the subset}%
\[
H(\mu ):=\dbigcup \mathcal{B}_{\Gamma }(\mu )=\dbigcup \{\gamma B:B\in 
\mathcal{B}_{\mathcal{H}}(\mu ),\gamma \in \Gamma \}. 
\]

\noindent \textbf{Proof. }Take a countable subgroup $\Gamma $ in $G,\ $which
is dense in $G$ under $\mathcal{\tau }_{G}$ and contains $D(\mathcal{H})$,
as in Lemma 7.4. Consider $w\in \gamma B_{\delta }^{K,\Delta }$ with $\gamma
\in \Gamma ,$ $K\in \mathcal{H},\delta ,\Delta \in \mathbb{Q}_{+},\mathit{\ }%
\Delta <\mu (K).$ Then for some $g\in D(\mathcal{H})\cap B_{\delta
}^{K,\Delta }$ and $\varepsilon >0$%
\[
gB_{\varepsilon }^{Kg,\Delta }\subseteq B_{\delta }^{K,\Delta }. 
\]%
So both%
\[
w\in \gamma gB_{\varepsilon }^{Kg,\Delta }\subseteq \gamma B_{\delta
}^{K,\Delta }, 
\]%
and $\gamma g\in \Gamma ,$ the latter as $g\in D(\mathcal{H})\subseteq
\Gamma .$ So, by the Corollary 7.1 (of Lemma 7.2), the family $\mathcal{B}%
(\mu )$ forms a sub-base for a topology on the set of points%
\[
\dbigcup \{\gamma B:B\in \mathcal{B}(\mu ),g\in \Gamma \}.\qquad \square 
\]

\noindent \textbf{Remark. }The same proof shows that one may drop
countability in the conditions and second-countability in the conclusions.

\bigskip

\noindent \textbf{Definition. }We term the second-countable topology of the
preceeding theorem the $\mu $\textit{-topology.}

\bigskip

In the next result we take $\mu =\sigma (\mathbf{t}).$ As $1_{G}\in
B_{\delta }^{K,\Delta }\subseteq B_{\delta },$ the $\sigma (\mathbf{t})$%
-topology evidently refines the \textit{original topology} $\tau _{G}$ of $%
G. $ The finer topology could be discrete; in cases of interest, however,
this will not happen:

\bigskip

\noindent \textbf{Proposition 7.1 (Refinement). }\textit{For }$G$\textit{\
amenable at }$1$\textit{,\ }$\mathbf{t}$\textbf{\ }\textit{null and
non-trivial, the open sets }$B_{\delta }^{K,\Delta }$ \textit{of the }$%
\sigma (\mathbf{t})$\textit{-topology are infinite and refine the topology
induced by }$\tau _{G}$ \textit{on }$H(\sigma (\mathbf{t})).$

\bigskip

\noindent \textbf{Proof.} For $\{g_{n}\}_{n\in \mathbb{N}}$ dense in $G,$
write $D:=\{g_{n}:n\in \mathbb{N}\}.$ The open sets of $G$ are generated as
unions of sets of the form $gV,$ with $g\in D$ and $V$ an open nhd of $%
1_{G}. $ We show that these sub-basic sets of the $\sigma (\mathbf{t})$%
-topology refine the $\tau _{G}$-nhds of the identity. For $V$ a non-empty $%
\tau _{G}$-open nhd of $1_{G}$ choose $U$ to be a non-empty $\tau _{G}$-open
nhd of $1_{G}$, with $U^{-1}U\subseteq V$.

Consider any non-trivial null sequence $\mathbf{t}$ and, referring to the
Subcontinuity theorem (Th. 7.1), consider $\sigma =\sigma (\mathbf{t})\in 
\mathcal{P}(G).$ For $\{g_{n}\}$ dense in $G$, there is $n$ with $\sigma
(g_{n}U)>0;$ for otherwise, $\sigma (g_{n}U)=0$ for each $n$ and, since $%
G=\dbigcup\nolimits_{n}g_{n}U,$ we reach the contradiction $\sigma (G)=0$.
Pick $n$ with $\sigma (g_{n}U)>0;$ write $g$ for $g_{n}.$

Choose compact sets $K_{n}$ such that $\sigma $ is concentrated on $%
\dbigcup\nolimits_{n}K_{n}$ and a countable base $\mathcal{B}$ for $\tau
_{G}.$ Since%
\[
\sigma (gU)=\sigma (\dbigcup\nolimits_{m}K_{m}\cap gU)=\sigma
(\dbigcup\nolimits_{m}\{K_{m}\cap g\bar{B}:\bar{B}\subseteq U,B\in \mathcal{B%
},m\in \mathbb{N}\}), 
\]%
there are $m\in \mathbb{N}$ and $B\in \mathcal{B}$ with\textit{\ }$\mu
(K_{m}\cap g\bar{B})>0.$

Take $K:=K_{m}\cap g\bar{B};$\textit{\ }there is a subsequence\textit{\ }$%
\mathbf{s}=\mathbf{s}(K):=\{t_{m(k)}\}$ with%
\[
\sigma (Kt_{m(k)})>\sigma (K)/2\qquad \text{(}k\in \mathbb{N}\text{)},\text{
so}\qquad \mu _{\_}^{\mathbf{s}}(K)>0. 
\]%
So as $\Delta :=\sigma _{\_}^{\mathbf{s}}(K)/4<\sigma _{\_}^{\mathbf{s}%
}(K)<\sigma (K),$ by [BinO7, Lemma 1], there is $\delta >0$ with%
\[
1_{G}\in B_{\delta }^{K,\Delta }\subseteq K^{-1}K\subseteq \bar{B}%
^{-1}g^{-1}g\bar{B}\subseteq U^{-1}U\subseteq V, 
\]%
and $t_{m(n)}\in B_{\delta }^{K,\Delta }$ for all large enough $n$ as in
[BinO7, Lemma 1]. $\square $

\bigskip

\noindent \textbf{Remark. }Proposition 7.1 is connected to the
Steinhaus-Weil Theorem, Theorem SW, in \S 6\ above: a similar argument
gives, for $E$ non-left-Haar-null, that $1_{G}\in \mathrm{int}_{G}(\hat{E})$
for%
\[
\hat{E}:=\bigcup\nolimits_{\delta ,\Delta >0,g\in G,\mathbf{t}}\{B_{\delta
}^{gK,\Delta }(\sigma (\mathbf{t})):K\subseteq E,K\in \mathcal{K}_{+}(\sigma
(\mathbf{t})),\Delta <\sigma (gK)\}. 
\]%
That is, the relevant basic open nhds of $1_{G}$ in the various $\sigma (%
\mathbf{t})$-topologies `aggregate' to yield a nhd of $1_{G}$ in the
original topology of $G.$

\section{Connections with Cameron-Martin theory}

In this section, we pursue the connection with Cameron-Martin theory.
Proposition 8.1 provides the basis for a definition of the `covered points'
under $\mu ;$ this identifies a canonical `largest' Steinhaus support for $%
\mu $ (modulo an initial choice of dense subset). The result takes its
motivation from the following classical observation:

\textit{In a locally convex topological vector space }$X$\textit{, in
particular in a Fr\'{e}chet space, equipped with a symmetric Radon Gaussian
measure }$\gamma $\textit{: if }$E$\textit{\ is any Hilbert space
continuously embedded in }$H(\gamma )$\textit{, then there exists a
symmetric Radon Gaussian measure }$\gamma ^{\prime }$\textit{\ with }$%
H(\gamma ^{\prime })=E$ [Bog1, 3.3.5].

We sign off by showing that the topology of the Steinhaus support is
metrizable.

We recall from \S 4.2: $\mu _{g}(B)=\mu (Bg)$ for $B\in \mathcal{B}(G);$%
\[
\mathcal{N}\mathcal{=}\mathcal{N}(\mu ):=\{g\in G:\mu _{g}\bot \mu \},\qquad
G_{0}=G_{0}(D):=\dbigcup\nolimits_{d\in D}d\mathcal{N}(\mu ), 
\]%
with $D=\{g_{n}:n\in \mathbb{N}\}\ $a dense subset of $G$; measures of
Solecki type (\S 7) have $\Delta (K,\delta )\equiv 0.$

\bigskip

\noindent \textbf{Proposition 8.1 (Covering Lemma).} \textit{For }$\mu \in 
\mathcal{P}(G)$ \textit{of Solecki type, let }$\tilde{D}$\textit{\ be a
dense subset of }$Q(\mu ).$\textit{\ For }$\delta >0,K\in \mathcal{K}%
_{+}(\mu )$ \textit{and }$x\in Q(\mu )$\textit{\ there is }$g\in \tilde{D}$%
\textit{\ with}%
\[
x\in gB_{\delta }^{K,\Delta } 
\]%
\textit{\ for all small enough }$\Delta <\mu (K)$\textit{.}

\bigskip

\noindent \textbf{Proof. }Choose $g\in \tilde{D}\cap xB_{\delta }\subseteq
Q(\mu )$. Then $y^{-1}:=x^{-1}g\in B_{\delta },$ so also $y=g^{-1}x\in
B_{\delta }$ (symmetry of the group-norm on $G),$ and $y=g^{-1}x\in Q(\mu )$%
, as $Q(\mu )$ is a subgroup. Now $\mu (Ky)>0,$ as $y\in Q(\mu ),$ so we may
choose $0<\Delta <\min \{\mu (Ky),\mu (K));$ then%
\[
x=gy\in g\{z\in B_{\delta }:\mu (Kz)>\Delta \}=gB_{\delta }^{K,\Delta
}.\qquad \square 
\]

Proposition 8.1 above identifies how points of $Q(\mu )$ can be covered by
certain translates of basic sets of the form $B_{\delta }^{K,\Delta }.$ To
go beyond $Q(\mu )$ this motivates the following.

\bigskip

\noindent \textbf{Definition.} Say that $g\in G$ is a \textit{covered} point
($g$ `is covered') under $\mu \in \mathcal{P}(G)$ if there is $K\in \mathcal{%
K}_{+}(\mu )$ with $\mu (Kg)>0.$ (Then $g\in B_{\delta }^{K,\Delta }$ for $%
\delta >||g||$ and $0<\Delta <\min \{\mu (K),\mu (Kg)\}.)$ So the points of $%
Q_{R}$ are covered, but $g\in G$ is \textit{not covered} if $\mu _{g}(K)=\mu
(Kg)=0$ for all $K\in \mathcal{K}_{+}(\mu ),$ that is, $\mu _{g}\bot \mu $.

\bigskip

\noindent \textbf{Proposition 8.2. }\textit{For }$\mu \in \mathcal{P}(G)$ 
\textit{of Solecki type, }$\{g_{n}\}_{n\in \mathbb{N}}$\textit{\ dense in }$%
G $ \textit{and }$\delta >0,$\textit{\ the sets }$g_{n}B_{\delta }^{K,\Delta
}$\textit{\ cover }$G\backslash G_{0}=G\backslash \dbigcup\nolimits_{n\in 
\mathbb{N}}g_{n}\mathcal{N}(\mu )$ \textit{and so generate a topology on the
Borel set }$G\backslash G_{0}$\textit{\ for which these are sub-basic.}

\bigskip

\noindent \textbf{Proof. }Consider $x\notin \dbigcup\nolimits_{n\in \mathbb{N%
}}g_{n}\mathcal{N}(\mu )$ with $\{g_{n}\}$ is dense in $G$. Then for
arbitrary $\delta >0$, select $g_{n}\in B_{\delta }(x).$ Then $x\in
B_{\delta }(g_{n}),$ and $y:=g_{n}^{-1}x\in B_{\delta }\backslash \mathcal{N}%
(\mu ),$ as $g_{n}^{-1}x\notin \mathcal{N}(\mu ).$ Now we may choose $K\in 
\mathcal{K}_{+}(\mu )$ with $\mu (Ky)>0.$ Then for $0<\Delta <\min \{\mu
(K),\mu (Ky)\},$%
\[
x=g_{n}y\in g_{n}B_{\delta }^{K,\Delta }. 
\]%
That is, for $\delta >0,$ the family $\{g_{n}B_{\delta }^{K,\Delta }:K\in 
\mathcal{K}_{+}(\mu ),$ $0<\Delta <\mu (K),n\in \mathbb{N}\}$ covers $%
G\backslash G_{0},$ and so a second-countable topology is generated with
sub-base the sets%
\[
\{g_{n}B_{\delta }^{K,\Delta }:K\in \mathcal{K}_{+}(\mu ),0<\Delta <\mu
(K),n\in \mathbb{N},\delta \in \mathbb{Q}_{+}\}.\qquad \square 
\]

\bigskip

\noindent \textbf{Remark. }In the special case when $Q(\mu )$ is dense in $G$
(for instance taking $G$ to be $\bar{Q}),$ so that also $\tilde{D}$ (in
Prop. 8.1) is dense in $G,$ Prop. 8.2 above (with $g_{n}=\tilde{g}_{n})$
follows from Proposition 8.1. Note that in this case also $G_{0}=\mathcal{N}%
(\mu ),$ since $\tilde{g}_{n}\mathcal{N}(\mu )\subseteq \mathcal{N}(\mu ),$
by Lemma 3.2.

\bigskip

\noindent \textbf{Definitions.} For $\mu \in \mathcal{P}(G)$ and $K\in 
\mathcal{K}_{+}(\mu ),$ put%
\[
\Delta _{K}(x,y):=|\mu (Kx)-\mu (Ky)|\leq 1\quad \quad (x,y\in G), 
\]%
which is a pseudometric, so that%
\[
\rho _{K}(x,y):=\max \{d_{L}^{G}(x,y),\Delta _{K}(x,y)\}\quad \quad (x,y\in
G) 
\]%
is a metric.

\bigskip

\noindent \textbf{Proposition 8.3. }\textit{For }$\mu \in \mathcal{P}(G)$ 
\textit{of Solecki type,} $g\in G\backslash \mathcal{N}(\mu )$, $K\in 
\mathcal{K}_{+}(\mu )$ \textit{and }$\varepsilon >0:$\textit{\ if }$%
0<\varepsilon <\mu (Kg),$\textit{\ then there is }$\delta =\delta
(\varepsilon )$\textit{\ with }$0<\delta <\varepsilon $\textit{\ such that}%
\textbf{\ }%
\[
gB_{\delta (\varepsilon )}^{Kg,\Delta }\subseteq B_{\varepsilon }^{\rho
_{K}}(g):=\{x:\rho _{K}(x,g)<\varepsilon \}\subseteq gB_{\varepsilon
}^{Kg,\Delta }\text{, for }\Delta =\mu (Kg)-\varepsilon . 
\]%
\textit{Hence, for any enumeration }$\{K_{n}\}_{n}$ \textit{of the basic
compact sets in }$\mathcal{K}_{+}(\mu )$ \textit{comprising the family }$%
\mathcal{H}_{0}$ \textit{of Lemma 7.4, the metric}%
\[
\rho (x,y):=\sup \{d_{L}^{G}(x,y),2^{-n}\Delta _{K_{n}}(x,y)\} 
\]%
\textit{generates the }$\mu $\textit{-topology on }$G\backslash G_{0}$%
\textit{.}

\bigskip

\noindent \textbf{Proof. }Note that for $0<\varepsilon <\mu (Kg)$%
\[
\Delta _{K}(x,g)<\varepsilon \text{ iff [}d_{L}^{G}(x,g)<\varepsilon \text{
and }\mu (Kg)-\varepsilon <\mu (Kx)<\mu (Kg)+\varepsilon ]. 
\]%
Write $x=gh;$ then $d_{L}^{G}(x,g)<\varepsilon $ is equivalent to the
constraint $||h||<\varepsilon .$ As $x\rightarrow \mu (Kx)$ is upper
semicontinuous, there is $0<\delta =\delta (\varepsilon )<\varepsilon $ such
that $\mu (Kgh)<\mu (Kg)+\varepsilon ,$ for $h\in B_{\delta };$ this yields
the further required constraint $\mu (Kgh)>\Delta :=\mu (Kg)-\varepsilon .$
The remaining assertions are now immediate. $\square $

\bigskip

\noindent \textbf{Remarks.} 1. In the above argument $\mu (Kgh\triangle
Kg)\leq 2\varepsilon $ provided $\mu (KgB_{\delta })<\varepsilon .$ This
implies that convergence in $\rho $ implies convergence in the Weil-like
group-norm $||\cdot ||_{\mu }^{E}$ of [BinO7] with $E=Kg$; indeed in the
locally compact case these norms generate the Weil topology of [Wei] (cf.
[Hal, \S 62], [HewR, \S 16], [Yam2, Ch. 2] and the recent [BinO7]). So the $%
\rho $-topology refines the Weil-like topology.

\noindent 2. As with Theorem 8.1 above, there is an analogue, in which
metrizability is dropped in favour of a uniform structure.

\section{Complements}

\noindent \textbf{1.} \textit{Historical remarks}. The fundamental reference
here is of course the first, Haar's 1933 paper in which he introduces Haar
measure [Haa]. Von Neumann, in the paper (of the same journal) immediately
after Haar's [Neu1], applies Haar measure for compact Lie groups to solve
Hilbert's fifth problem. He follows this with two further contributions
[Neu3,4]. Kakutani made extensive relevant contributions to both topological
groups and to measure theory. His papers on the first appear in Volume 1 of
his Selected Papers [Kak3], together with commentaries (p. 391-408) by A. H.
Stone, J. R. Choksi, W. W. Comfort, K. A. Ross and J. Taylor. Here he deals
with metrisation, with uniqueness of Haar measure, and (with Kodaira) on its
completion regularity. His papers on the second appear in Volume 2, with
commentaries (p.379-383) by Choksi, M. M. Rao, Oxtoby [Oxt3], and Ross. Here
he deals (alone, with Kodaira, and with Oxtoby) on extension of Lebesgue
measure, and with equivalence of infinite product measures (\S 9.18 below).

The other key historical references here are the Cameron-Martin papers
[CamM1,2,3]; see [Bog1], [LedT], [Str] for textbook accounts.

\noindent \textbf{2. }\textit{Radon measures. }We recall that on a complete
separable metric space every Borel measure is Radon, i.e. has inner compact
regularity ([Bog2, Vol. II, Th. 7.1.7]).

A metrizable \v{C}ech-complete space (i.e. one that is a $\mathcal{G}%
_{\delta }$ in some (any) compactification) is topologically complete, i.e.
the topology may be generated by a complete metric [Eng, Th. 4.3.26]. So, in
particular, if a locally compact group is metrizable, then it has a complete
metric, and so every Borel measure on the group, in particular every Haar
measure, is Radon. If in addition the group is separable (so Polish), then,
being second-countable, it is $\sigma $-compact, and then every Haar measure
is $\sigma $-finite, and so also outer regular ([Kal, Lemma 1.34], cf. [Par,
Th. II.1.2] albeit for a probability measure).

In general, one may pass from a Haar measure $\eta _{X}$ on a locally
compact group $X$ which is outer regular (i.e. Borel sets are outer $\eta
_{X}$-approximable by open sets) to the Borel measure $\mu $ defined by%
\[
\mu (B):=\sup \{\eta _{X}(K):K\in \mathcal{K}(X),K\subseteq B\}\quad (B\in 
\mathcal{B}(X)), 
\]%
which agrees with $\eta _{X}$ on $\mathcal{K}(X)$ and so is inner compact
regular [Bog2, Th. 7.11.1]; however, $\mu $ need not be outer-regular. In
applications inner compact regularity carries more advantages, hence our
adoption of this property of measures.

We note some alternative usages here.

(a) For Schwartz [Sch, 1.2], a Radon measure is a locally finite, Borel
measure which is inner compact regular (definition R$_{\text{3}}$); an
equivalent definition includes local finiteness and couples outer regularity
with inner compact regularity restricted to open sets (definition R$_{\text{2%
}}$). A third definition involves both a locally finite measure $M$ which is
outer regular and $m$ its associated essential measure (outer measure
restricted to $\mathcal{B}(X)$) which is inner compact regular, the two
agreeing on open sets and on sets of finite $M$-measure.

(b) Fremlin [Fre2, p. 15] defines Radon measures to be locally finite and
inner compact regular (plus complete and locally determined [Fre1, p.13]).

(c) Heyer [Hey1, \S 1.1], for a locally compact group $X$, defines a Radon
measure as a linear functional with domain the continuous complex functions
with compact support in $X$ and with a boundedness condition where the
bounds correspond to the possible compact supports.

\noindent \textbf{3. }\textit{Invariance beyond local compactness. }We
recall our opening paragraph, which set out the contrast between the local
compactness of the group setting, where one has Haar measure, and the
absence of both in the Hilbert-space setting in which Cameron-Martin theory
originated. We note that the invariance property of Haar measure may be
extended \textit{beyond }the locally-compact case. Nothing new is obtained
in our setting of probability measures, but if one drops local finiteness,
Haar-like measures of `pathological' character can occur (\S 9.8 below). We
quote Diestel and Spalsbury [DieS, Ch. 10], who give a textbook account of
the early work of Oxtoby in this area [Oxt1]. We note in passing that this
interesting paper is not cited by Oxtoby himself in either edition of his
classic book [Oxt2]. We note also the use of local finiteness in Schwartz's
definition of a Radon measure [Sch].

The classic case of Haar (invariant) measure is Lebesgue measure in
Euclidean space. A number of authors have produced `Lebesgue-like'
extensions of Lebesgue measure from $\mathbb{R}^{n}$ to $\mathbb{R}^{\mathbb{%
N}};$ see e.g. Baker [Bak1,2], Gill and Zachary [GilZ], [Pan], Yamasaki
[Yam1,2].

Admissible translators present themselves here and also in a variety of
related circumstances; for a statistical setting see e.g. [Shep], and for
later developments [Smol] and Sadasue [Sad1].

\noindent \textbf{4. }\textit{Quasi-invariance beyond local compactness. }%
Such questions are addressed in a vector-space context in Bogachev's book
[Bog1]; see also Yamasaki [Yam1,2], Arutyunyan and Kosov [AruK] (cf. \S 4).
For the group context, see Ludkovsky [Lud], Sadasue [Sad1,2] and the
references cited there (again, cf. \S 4).

\noindent \textbf{5. }\textit{Group representations beyond local
compactness. }See Ludkovsky [Lud] for group representations, and the
monograph of Banaszczyk [Bana] for Pontryagin duality in the abelian case;
for a textbook treatment see [FelD]. For harmonic analysis, see Gel$^{\prime
}$fand and Vilenkin [GelV], Xia [Xia].

\noindent \textbf{6. }\textit{Integration beyond local compactness. }Measure
and integration are of course closely linked, in this context as in any
other. For monograph accounts, see e.g. Skorohod [Sko], [Yam2], [Xia],
[GilZ].

\noindent \textbf{7. }\textit{Differentiation beyond local compactness. }%
Differentiation in infinitely many dimensions owes much to pioneering work
by Fomin, and has led on to the theory of smooth measures and the Malliavin
calculus. See e.g. Bogachev [Bog3], Dalecky and Fomin [DalF].

\noindent \textbf{8. }\textit{The Oxtoby-Ulam Theorem }([Oxt1, Th. 2],
[DieS, Th.10.1]). This asserts that in a non-locally-compact Polish group
carrying a (non-trivial, left) invariant Borel measure every nhd of the
identity contains \textit{uncountably} many \textit{disjoint} (left) \textit{%
translates} of a compact set of positive measure. Since local finiteness
rules out such pathology, `total' invariance of a Radon measure implies
local compactness, hence the introduction of `selective invariance' and
`selective approximation' (by compact sets) in the variant Steinhaus triples
of \S 2.

\noindent \textbf{9. }\textit{Invariant means. }One can deal with invariant
means in place of invariant measures. This involves the theory of \textit{%
amenable groups, }and amenability more generally; see Paterson [Pat], which
has an extensive bibliography. There are links with Solecki's concept of 
\textit{amenability at }1 ([Sol] and \S 6; [BinO7]).

\noindent \textbf{10. }\textit{Fr\'{e}chet spaces: Gaussianity}\textbf{\ }%
[Bog1]. For $X$ a locally convex topological vector space, $\gamma $ a
probability measure on the $\sigma $-algebra of the cylinder sets generated
by $X^{\ast }$ (the Borel sets, for $X$ separable Fr\'{e}chet, e.g.
separable Banach), with $X^{\ast }\subseteq L^{2}(\gamma ):$ then $\gamma $
is called \textit{Gaussian} on $X$ iff $\gamma \circ \ell ^{-1}$ defined by 
\[
\gamma \circ \ell ^{-1}(B)=\gamma (\ell ^{-1}(B))\qquad (\text{Borel }%
B\subseteq \mathbb{R}) 
\]%
is Gaussian (normal) on $\mathbb{R}$ for every $\ell \in X^{\ast }\subseteq
L^{2}(\gamma )$. For a monograph treatment of Gaussianity in a Hilbert-space
setting, see Janson [Jan].

\noindent \textbf{11. }\textit{Cameron-Martin} \textit{aspects. }For $X\ $Fr%
\'{e}chet and $\gamma $ Gaussian case, when the closed span of $\{x^{\ast
}-\gamma (x^{\ast }):x^{\ast }\in X^{\ast }\}$ is infinite-dimensional, $%
H(\gamma )$ is $\gamma $-null in $X$ [Bog1, Th. 2.4.7].

Furthermore, for $h\in H(\gamma ),$ the Radon-Nikodym density $d\gamma
_{h}/d\gamma $ (which is explicitly given by the Cameron-Martin formula $%
(CM) $) as a function of $h$ is continuous on $H(\gamma )$ [Bog1, Cor.
2.4.3]. This implies, for $\mathbf{t}$ null in the $H(\gamma )$-norm (with $%
t_{n}\in H(\gamma ))$ and compact $K$ with $\gamma (K)>0,$ that $\gamma
_{-}^{\mathbf{t}}(K)>0.$ Here, for $X$ sequentially complete, the
corresponding balls (under the $H(\gamma )$-norm) are weakly compact in $X$
-- cf. [Bog1, Prop. 2.4.6], also the Remark before Th. 3.4 above.

We note here for convenience the following properties of the Cameron-Martin
space.

(i) For $\gamma $ non-degenerate, $H(\gamma )$ is everywhere dense [Bog1, \S %
3.6].

(ii) $\mathrm{cl}_{X}H(\gamma )$ is of full measure [Bog1, 3.6.1].

(iii) If $X_{\gamma }^{\ast }$ is infinite-dimensional (i.e. is not locally
compact), then $\gamma (H(\gamma ))=0$ [Bog1, 2.4.7].

(iv) For $\gamma $ a Radon Gaussian measure, both of the spaces $H(\gamma )$
and $L^{2}(X,\gamma )$ are separable [Bog1, 3.2.7 and Cor. 3.2.8].

(v) The `relative mobility property' (cf. \S 9.15 below) that $\gamma
(Kh)\rightarrow \gamma (K)$ as $h\rightarrow 0$ always holds [Bog1, Th.
2.4.8] applied to $1_{K}$ (indeed, for $h\in H(\gamma ),$ $d\gamma
_{h}/d\gamma $ exists and is continuous in $h$, by [Bog1, 2.4.3]).

(vi) For $\gamma $ Radon and $X$ a locally convex topological vector space,
the closed unit ball of $H(\gamma )$ is compact in $X$ [Bog1, 3.2.4].

(vii) In a locally convex space, there is a sequence of metrizable compacta $%
K_{n}$ with $\gamma (\tbigcup\nolimits_{n}K_{n})=1$ [Bog1, Th. 3.4.1].

(viii) For $X$ a locally convex space, equipped with a symmetric Radon
Gaussian measure $\gamma $: if $E$ is any Hilbert space continuously
embedded in $H(\gamma )$, then there exists a symmetric Radon Gaussian
measure $\gamma ^{\prime }$ with $H(\gamma ^{\prime })=E$ [Bog1, 3.3.5].

\noindent \textbf{12. }\textit{Locally compact groups: Gaussianity. }For
Gaussian measures on locally compact groups $G$, see e.g. [Par, IV.6] for $G$
abelian and [Hey1, 5.2, 5.3] for the general case. Use is made there of
characters -- bounded, multiplicative or additive according to notation; the
local inner product [Hey1, 5.1.7] is between the group and its Pontryagin
dual.

One link between the group and vector-space aspects can be seen in the
central role played in each by Gaussianity. We may think of this in each
case as saying that, as in $(CM),$ the relevant Fourier transform is of
exponential type, the exponent having two terms, one linear (concerning
means -- location, or translation), one quadratic (concerning covariances,
which captures scale and dependence effects). Where the density of the
measure exists, it involves (via the `normal' Edgeworth formula above in the
Euclidean case) the inverse of the covariance $\Sigma $ (matrix or
operator), important in its own right (as the concentration or precision
matrix/operator $K:=\Sigma ^{-1}$). So `degeneracy-support' phenomena as
above are unavoidable (below). Statistically, samples from two populations
can only be usefully compared if their covariances are the same, and then
the relevant statistic is the likelihood ratio; see e.g. [IbrR] for
background here.

The supports of Gaussian measures on groups, and in particular the
connections between Gaussian and Haar supports, have been studied by
McCrudden [McC1,2], [McCW] and others\textit{.}

\noindent \textbf{13. }\textit{Dichotomy. }The equivalence-singularity
dichotomy for Gaussian measures is a general consequence ([LePM], [Kal];
[MarR, \S 5.3]) of the triviality of a certain tail algebra (cf. [Hey1, Th.
5.5.6]); tail triviality in this case is established using a \textit{%
zero-one law}.

\noindent \textbf{14.}\textit{\ Automatic continuity. }The general theme of
automatic continuity -- situations where a function subject to mild
qualitative conditions must necessarily be continuous -- is important in
many contexts; see e.g. [BinO5] and the references therein. For results of
this type on $\gamma $-measurable linear functions for Gaussian $\gamma $,
see [Bog1, Ch. 2]. See also [Sol, Cor. 2].

\noindent \textbf{15.}\textit{\ Simmons-Mospan theorem and subcontinuity. }%
Recall from the companion paper [BinO7] (cf. [BinO4]) the following
definition, already used in Th. 3.4 above. For $\mu \in \mathcal{P}(G)$ and $%
K\in \mathcal{K}(G),$ 
\[
\mu _{-}(K):=\sup_{\delta >0}\inf \{\mu (Kt):t\in B_{\delta }\}; 
\]%
then $\mu $ is \textit{subcontinuous} if $\mu _{\_}(K)>0$ for all $K$ with $%
\mu (K)>0.$ (For a related notion see [LiuR], where a Radon measure $\mu $
on a space $X,$ on which a group $G\ $acts homeomorphically, is called 
\textit{mobile} if each map $t\mapsto \mu (Kt)$ is continuous for $K\in 
\mathcal{K}(X).)$ It follows from Prop. 3.1 above (see Cor. 3.2) that if $%
\mu _{-}(K)>0$ for some $K,$ then $G$ is locally compact. Note that in a
locally compact group, right uniform continuity of all the maps $t\mapsto
\mu (tB)$ for $B$ Borel is equivalent to absolute continuity of $\mu $
w.r.t. Haar measure ([Hey1, L. 6.3.4] -- cf. [HewR, Th. 20.4]). So if $G\ $%
is not locally compact, no measure $\mu \in \mathcal{P}(G)$ is
subcontinuous; then for all compact $K\subseteq G,$ $\mu _{-}(K)=0.$

If $G$ is locally compact, then its left Haar measure $\eta =\eta _{G}$
satisfies 
\[
\eta _{-}(K)=\eta (K); 
\]%
in particular, this equation holds for all non-null compact $K.$ The latter
observation extends to measures $\mu $ that are absolutely continuous w.r.t. 
$\eta _{G}.$ Conversely, if $\mu $ is a measure satisfying $\mu _{-}(K)>0$%
\textit{\ }for all compact $K$ with $\mu (K)>0,$ then, as a consequence of
the Simmons-Mospan theorem (\S 1), $\mu $ is absolutely continuous w.r.t. $%
\eta _{G}:$ see [BinO7].

\textbf{16. }\textit{Quasi-invariance and the Mackey topology of analytic
Borel groups. }We stop to comment on the force of \textit{full
quasi-invariance} of a measure in connection with a Steinhaus triple $%
(H,G,\mu )$ with $H$ (and $G$) Polish. Both groups, being absolutely Borel,
are analytic spaces (Lemma 2.1 above). So both carry a \textit{standard}
Borel `structure' (i.e. Borel isomorphic to the $\sigma $-algebra of Borel
subsets of some Borel subset of a Polish space) with $H$ carrying a Borel
`substructure' ($\sigma $-subalgebra) of $G.$ (Borel subsets of $H$ are
Borel in $G.$) Mackey [Mac] investigates such Borel groups, defining also a
(Borel) measure $\mu $ to be \textit{standard} if it has a \textit{standard}
Borel support (vanishes outside of a standard Borel set). It emerges that
every $\sigma $-finite Borel measure in an analytic Borel space is standard
[Mac, Th. 6.1]. Of interest to us is Mackey's notion of a `measure class' $%
C_{\mu },$ comprising all Borel measures $\nu $ with the same null sets as $%
\mu :$ $\mathcal{M}_{0}(\nu )=\mathcal{M}_{0}(\mu ).$ Such a measure class
may be closed under translation, and may be right or left invariant; then
the common null sets are themselves invariant, and so may be viewed as
witnessing quasi-invariance of the measure $\mu .$ Mackey shows that a Borel
group with a one-sided invariant measure class has a both-sidedly invariant
measure class [Mac, Lemma 7.2]; furthermore, if the class is countably
generated, then the class contains a left-invariant and a right-invariant
measure [Mac, Lemma 7.3]. This enables Mackey to improve on Weil's theorem
in showing that an analytic Borel group $G$ with a one-sidedly invariant
measure class, in particular one generated by a quasi-invariant measure, has
a unique locally compact topology making $G$ a topological group as well as
generating the given Borel `structure'.

\textbf{17. }\textit{The Strassen set and the law of the iterated logarithm
(LIL)}. The LIL completes (with the law of large numbers (LLN) and central
limit theorem (CLT)) the trilogy of classical limit theorems in probability
theory; for a survey see e.g. [Bin]. One form, the compact LIL, links the
unit ball $U$ of the reproducing-kernel Hilbert space associated with the
covariance structure of a random variable $X$ with values in a separable
Banach space $B$ with the cluster set of the partial sums, normalised as in
the classical LIL. See e.g. [LedT, 207-210]. The first results of this type
were Strassen's functional LIL and its extension to Banach spaces by Kuelbs
and others ([LedT, 233-234], [Bog1, 358]).

\textbf{18. }\textit{Product measures}. Infinite products of probability
measures correspond to infinite sequences of independent random variables;
they give a particularly important class of probability measures on
infinite-dimensional spaces. A basic result here is the \textit{Kakutani
alternative}: if the laws of the factors are equivalent, the laws of the
products are either equivalent or mutually singular, depending on the
convergence or divergence of the infinite product of the Hellinger distances
of the factor laws ([Kak2], [JacS, Ch. IV]; the term \textit{%
Kakutani-Hellinger distance} is now used). As usual with dichotomies in
probability theory, there are links with zero-one laws (cf. \S 9.13). See
also Shepp [She], [Kak2].

\textbf{19. }\textit{Non-Archimedean fields}. L\"{o}wner's result [Loe] (cf.
[Neu2]) addresses the loss of a property -- desirable, in some respects --
as the dimension $n\rightarrow \infty $, by changing the base field from the
reals to a non-Archimedean field. This is an early example of
non-Archimedean fields (which originate in algebra and algebraic number
theory) being applied to address a concrete problem in a quite different
area.

\textbf{20. }\textit{Other settings. }Recent generalizations of
Cameron-Martin theory analyze an infinite-dimensional Lie-group or a
sub-Riemannian manifold setting -- see for example [DriG], [GorL] and [Gor],
which thus preserve much of the classical setting; see also [Pug] (cf.
[Bog1, p. 393]) and [Shi] for special cases.

\bigskip

\textbf{Acknowledgement. }We thank V. I. Bogachev for helpful discussion and
the Referee for his thorough, scholarly and constructive report, which led
to many improvements.

\bigskip

\textbf{References}

\bigskip

\noindent \lbrack ArhT] A. Arhangel'skii, M. Tkachenko, \textsl{Topological
groups and related structures}. World Scientific, 2008. \newline
\noindent \lbrack AruK] L. M. Arutyunyan and E. D. Kosov, Spaces of
quasi-invariance of product measures. \textsl{Funct. Anal. Appl.}\textbf{\ 49%
} (2015), 142--144.\newline
\noindent \lbrack Bad] R. Badora, On the Hahn-Banach theorem for groups. 
\textsl{Arch. Math.} (Basel) \textbf{86} (2006), 517--528.\newline
\noindent \lbrack Ban1] S. Banach, \textsl{Th\'{e}orie des op\'{e}rations lin%
\'{e}aires}, in: Monografie Mat., vol.1, 1932 (in: \textquotedblleft
Oeuvres\textquotedblright , vol.2, PWN, 1979), translated as `Theory of
linear operations, North Holland, 1978.\newline
\noindent \lbrack Ban2] S. Banach, Sur l'\'{e}quation fonctionnelle $%
f(x+y)=f(x)+f(y)$, \textsl{Fund. Math. }\textbf{1} (1920),123--124;
reprinted in: Oeuvres, Vol. I, PWN, Warszawa, 1967, 47--48 (commentary by H.
Fast, p. 314).\newline
\noindent \lbrack Bana] W. Banaszczyk, \textsl{Additive subgroups of
topological vector spaces.} Lecture Notes in Mathematics, 1466. Springer,
1991.\newline
\noindent \lbrack Bak1] R. L. Baker,\textquotedblleft Lebesgue
measure\textquotedblright\ on $\mathbb{R}^{\infty }$. Proc. Amer. Math. Soc.
113 (1991), no. 4, 1023--1029.\newline
\noindent \lbrack Bak2] R. L. Baker, \textquotedblleft Lebesgue
measure\textquotedblright\ on $\mathbb{R}^{\infty }$. II. \textsl{Proc.
Amer. Math. Soc.} \textbf{132} (2004), 2577--2591.\newline
\noindent \lbrack BarFF] A. Bartoszewicz, M. Filipczak, T. Filipczak, On
supports of probability Bernoulli-like measures. \textsl{J. Math. Anal. Appl.%
} \textbf{462} (2018), 26--35.\newline
\noindent \lbrack BecK] H. Becker, and A. S. Kechris, \textsl{The
descriptive set theory of Polish group actions.} London Math. Soc. Lecture
Notes \textbf{232}. Cambridge University Press, 1996.\newline
\noindent \lbrack BerTA] A. Berlinet, C. Thomas-Agnan, \textsl{%
Reproducing-kernel Hilbert spaces in probability and statistics}. Kluwer,
2004.\newline
\noindent \lbrack Ber] E. Berz, Sublinear functions on $\mathbb{R},$ \textrm{%
\textsl{Aequat. Math.}} \textbf{12} (1975), 200-206.\newline
\noindent \lbrack Bin] N. H. Bingham, Variants on the law of the iterated
logarithm. \textsl{Bull. London Math. Soc.} \textbf{18} (1986), 433-467.%
\newline
\noindent \lbrack BinF] N. H. Bingham and J.M. Fry, \textsl{Regression:
Linear models in statistics}. Springer, 2010.\newline
\noindent \lbrack BinGT] N. H. Bingham, C. M. Goldie and J. L. Teugels, 
\textsl{Regular variation}, 2$^{\text{nd}}$ ed., Cambridge University Press,
1989 (1$^{\text{st}}$ ed. 1987). \newline
\noindent \lbrack BinK] N. H. Bingham and R. Kiesel, \textsl{Risk-neutral
valuation. Pricing and hedging of financial derivatives. }Springer, 2$^{%
\text{nd}}$ ed. 2004 (1$^{\text{st}}$ed. 1998) \newline
\noindent \lbrack BinO1] N. H. Bingham and A. J. Ostaszewski, Kingman,
category and combinatorics. \textsl{Probability and mathematical genetics}
(Sir John Kingman Festschrift, ed. N. H. Bingham and C. M. Goldie), 135-168, 
\textsl{London Math. Soc. Lecture Notes in Mathematics} \textbf{378}, CUP,
2010. \newline
\noindent \lbrack BinO2] N. H. Bingham and A. J. Ostaszewski, Normed groups:
Dichotomy and duality. \textsl{Dissert. Math.} \textbf{472} (2010), 138p. 
\newline
\noindent \lbrack BinO3] N. H. Bingham, A. J. Ostaszewski, Category-measure
duality: convexity, midpoint convexity and Berz sublinearity. \textsl{%
Aequat. Math.} \textbf{91} (2017), 801--836.\newline
\noindent \lbrack BinO4] N. H. Bingham and A. J. Ostaszewski, Beyond
Lebesgue and Baire IV: Density topologies and a converse Steinhaus-Weil
theorem. \textsl{Topology and its Applications} \textbf{239} (2018), 274-292
(arXiv:1607.00031).\newline
\noindent \lbrack BinO5] N. H. Bingham and A. J. Ostaszewski, Additivity,
subadditivity and linearity: Automatic continuity and quantifier weakening. 
\textsl{Indag. Math.} (N.S.) \textbf{29} (2018), 687--713.\newline
\noindent \lbrack BinO6] N. H. Bingham and A. J. Ostaszewski, Variants on
the Berz sublinearity theorem, arXiv:1712.05183.\newline
\noindent \lbrack BinO7] N. H. Bingham and A. J. Ostaszewski, The
Steinhaus-Weil property: its converse, Solecki amenability and
subcontinuity, arXiv:1607.00049v3.\newline
\noindent \lbrack Bir] G. Birkhoff, A note on topological groups. \textsl{%
Compos. Math.} \textbf{3} (1936), 427--430.\newline
\noindent \lbrack Bog1] V. I. Bogachev, \textsl{Gaussian measures}, Math.
Surveys \& Monographs \textbf{62}, Amer Math Soc., 1998.\newline
\noindent \lbrack Bog2] V. I. Bogachev, \textsl{Measure theory.} Vol. I, II.
Springer, 2007.\newline
\noindent \lbrack Bog3] V. I. Bogachev, \textsl{Differentiable measures and
the Malliavin calculus.} Math. Surveys \& Monographs \textbf{164}. Amer Math
Soc., 2010.\newline
\noindent \lbrack CamM1] R. H. Cameron, W. T. Martin, Transformations of
Wiener integrals under translations. \textsl{Ann. Math. }(2) \textbf{45},
(1944), 386--396.\newline
\noindent \lbrack CamM2] R. H. Cameron, W. T. Martin, Transformations of
Wiener integrals under a general class of linear transformations. \textsl{%
Trans. Amer. Math. Soc.} \textbf{58}, (1945), 184--219.\newline
\noindent \lbrack CamM3] R. H. Cameron, W. T. Martin, The transformation of
Wiener integrals by nonlinear transformations. \textsl{Trans. Amer. Math.
Soc.} \textbf{66}, (1949), 253--283.\newline
\noindent \lbrack Chr1] J. P. R. Christensen, On sets of Haar measure zero
in abelian Polish groups. Proceedings of the International Symposium on
Partial Differential Equations and the Geometry of Normed Linear Spaces
(Jerusalem, 1972). \textsl{Israel J. Math.} \textbf{13} (1972), 255--260
(1973).\newline
\noindent \lbrack Chr2] J. P. R. Christensen, \textsl{Topology and Borel
structure. Descriptive topology and set theory with applications to
functional analysis and measure theory.} North-Holland Mathematics Studies 
\textbf{10}, 1974.\newline
\noindent \lbrack Con] J. B. Conway, \textsl{A course in functional analysis.%
} 2$^{\text{nd}}$ ed. Graduate Texts in Mathematics \textbf{96}. Springer,
1990 (1$^{\text{st}}$ ed. 1985).\newline
\noindent \lbrack DalF] Yu. L. Dalecky and S. V. Fomin, \textsl{Measures and
differential equations in infinite-dimensional space. }Kluwer, 1991.\newline
\noindent \lbrack DieS] J. Diestel, A. Spalsbury,\textsl{\ The joys of Haar
measure.} Graduate Studies in Mathematics \textbf{150}, Amer. Math. Soc.,
2014.\newline
\noindent \lbrack Dod] P. Dodos, The Steinhaus property and Haar-null sets. 
\textsl{Bull. Lond. Math. Soc.} \textbf{41} (2009), 377--384.\newline
\noindent \lbrack DriG] B. K. Driver, M. Gordina, Heat kernel analysis on
infinite-dimensional Heisenberg groups. \textsl{J. Funct. Anal.} \textbf{255}
(2008), 2395--2461.\newline
\noindent \lbrack FelD] J. M. G.Fell and R. S. Doran, \textsl{%
Representations of *-algebras, locally compact groups, and Banach
*-algebraic bundles}:\textsl{\ Basic representation theory of groups and
algebras.} Vol. 1 Pure and Applied Mathematics \textbf{125}. Academic Press,
1988.\newline
\noindent \lbrack For] {M. K. Fort, Jr., {A unified theory of
semi-continuity.} \textsl{Duke Math. J. }\textbf{16} (1949), 237--246.}%
\newline
\noindent \lbrack Fre1] D. H. Fremlin, \textsl{Measure theory}, Volume 2,
Broad foundations, Torres Fremlin, 2001.\newline
\noindent \lbrack Fre2] D. H. Fremlin, \textsl{Measure theory}, Volume 4,
Topological measure spaces, Part One, Torres Fremlin, 2003.\newline
\noindent \lbrack Fuc] L. Fuchs, \textsl{Infinite abelian groups.Vol. 1}
Academic Press, 1970; vol. 2 1973.\newline
\noindent \lbrack GarP] R. J. Gardner and W. F. Pfeffer, \textsl{Borel
measures}, in: \textsl{Handbook of set-theoretic topology} (K. Kunen and J.
E. Vaughan, eds), 961--1043, North-Holland, 1984.

\noindent \lbrack GelV] I. M. Gel$^{\prime }$fand and N. Ya. Vilenkin, 
\textsl{Generalized functions}. Vol. 4.\textsl{\ Applications of harmonic
analysis.} Academic Press , 1964.\newline
\noindent \lbrack GikS] I. I. Gikhman and A. V. Skorokhod, \textsl{Theory of
Random Processes I}, Grundlehren math. Wiss. \textbf{210}. Springer, 1974
(reprinted 2004; Izdat. Nauka, Moscow, 1971).\newline
\noindent \lbrack GilZ] T. Gill, W. Zachary, \textsl{Functional Analysis and
the Feynman Operator Calculus}, Springer, 2016.\newline
\noindent \lbrack Gir] I. V. Girsanov, On transforming a certain class of
stochastic processes by absolutely continuous substitution of measures, 
\textsl{Th. Prob. Appl}. \textbf{5} (1960), 285--301.\newline
\noindent \lbrack Gor] M. Gordina, Heat kernel analysis and Cameron-Martin
subgroup for infinite dimensional groups. \textsl{J. Funct. Anal.} \textbf{%
171} (2000), 192--232.\newline
\noindent \lbrack GorL] M. Gordina, T. Laetsch, A convergence to Brownian
motion on sub-Riemannian manifolds. \textsl{Trans. Amer. Math. Soc}. \textbf{%
369} (2017), 6263--6278.\newline
\noindent \lbrack Gow1] C. Gowrisankaran, Radon measures on groups. \textsl{%
Proc. Amer. Math. Soc.} \textbf{25} (1970), 381--384.\newline
\noindent \lbrack Gow2] C. Gowrisankaran, Quasi-invariant Radon measures on
groups. \textsl{Proc. Amer. Math. Soc. }\textbf{35} (1972), 503--506.\newline
\noindent \lbrack Gow3] C. Gowrisankaran, Semigroups with invariant Radon
measures. \textsl{Proc. Amer. Math. Soc.} \textbf{38} (1973), 400--404.%
\newline
\noindent \lbrack Gro1] L. Gross, Abstract Wiener spaces, in: \textsl{Proc.
Fifth Berkeley Sympos. Math. Statist. and Probability}, Vol. II, Part 1, pp.
31--42, Univ. California Press, 1967.\newline
\noindent \lbrack Gro2] L. Gross, Abstract Wiener measure and
infinite-dimensional potential theory, in: \textsl{Lectures in Modern
Analysis and Applications}, II, pp. 84--116, \textsl{Lecture Notes in
Mathematics} \textbf{140}, Springer, 1970.\newline
\noindent \lbrack Haa] A. Haar, Der Massbegriff in der Theorie der
kontinuierlichen Gruppen. \textsl{Math. Ann.} \textbf{34} (1933), 147-169.%
\newline
\noindent \lbrack Hal] P. R. Halmos{\normalsize , \textsl{Measure theory},}
Grad. Texts in Math. \textbf{18}, Springer 1974 (1st. ed. Van Nostrand,
1950).{\normalsize \newline
}\noindent \lbrack Hey1] H. Heyer, \textsl{Probability measures on locally
compact groups.} Ergebnisse Math. \textbf{94}, Springer, 1977.\newline
\noindent \lbrack Hey2] H. Heyer, Recent contributions to the embedding
problem for probability measures on a locally compact group. \textsl{J.
Multivariate Anal.} \textbf{19} (1986), 119--131.\newline
\noindent \lbrack HeyP] H. Heyer, G. Pap, On infinite divisibility and
embedding of probability measures on a locally compact abelian group. 
\textsl{Infinite-dimensional harmonic analysis} III, 99--118, World Sci.
Publ., 2005.\newline
\noindent \lbrack HewR] E. Hewitt and K. A Ross, \textsl{Abstract harmonic
analysis.} Vol. I., Grundlehren Math. Wiss. \textbf{115}, Springer, 2$^{%
\text{nd}}$ ed. 1979 (1$^{\text{st}}$ ed. 1963).\newline
\noindent \lbrack IbrR] I. A. Ibragimov, Y. A. Rozanov, \textsl{Gaussian
random processes.} Translated from the Russian by A. B. Aries. Applications
of Mathematics \textbf{9}. Springer, 1978.\newline
\noindent \lbrack JacS] J. Jacod and A. N. Shiryaev, \textsl{Limit theorems
for stochastic processes.}Grund. Math. Wiss. \textbf{288}, 2$^{\text{nd}}$
ed. 2003 (1$^{\text{st}}$ ed. 1987).\newline
\noindent \lbrack Jan] S. Janson, \textsl{Gaussian Hilbert spaces.}
Cambridge Tracts in Mathematics \textbf{129}. Cambridge University Press,
1997.\newline
\noindent \lbrack Kak1] S. Kakutani, \"{U}ber die Metrisation der
topologischen Gruppen, \textsl{Proc. Imp. Acad. Tokyo} \textbf{12} (1936),
82--84 (reprinted in [Kak3], Vol. 1, 60--62).\newline
\noindent \lbrack Kak2] S. Kakutani, On equivalence of infinite product
measures. \textsl{Ann. of Math.} (2) \textbf{49}, (1948). 214--224
(reprinted in [Kak3], Vol. 2, 19-29).\newline
\noindent \lbrack Kak3] S. Kakutani, \textsl{Selected Papers}, Vol. 1, 2,
ed. R. R. Kallman, Birkh\"{a}user, 1986.\newline
\noindent \lbrack Kal] G. Kallianpur, Zero-one laws for Gaussian processes. 
\textsl{Trans. Amer. Math. Soc.} \textbf{149} (1970), 199-211.\newline
\noindent \lbrack Kap] I. Kaplansky, \textsl{Infinite abelian groups.} U.
Michigan Press, 1969 (1$^{\text{st}}$ ed. 1954).\newline
\noindent \lbrack Kec] A. S. Kechris, \textsl{Classical descriptive set
theory}. Graduate Texts in Mathematics \textbf{156}, Springer, 1995.\newline
\noindent \lbrack Kem] J. H. B. Kemperman, A general functional equation.%
\textsl{\ Trans. Amer. Math. Soc.} \textbf{86} (1957), 28--56.{\normalsize 
\newline
}\noindent \lbrack Kle] V. Klee, Invariant extensions of linear functionals, 
\textsl{Pacific J. Math.} 4 (1954), 37-46.\newline
\noindent \lbrack Kuc] M. Kuczma, \textsl{An introduction to the theory of
functional equations and inequalities. Cauchy's equation and Jensen's
inequality,} 2$^{\text{nd}}$ ed., Birkh\"{a}user, 2009 (1$^{\text{st}}$ ed.
PWN, Warszawa, 1985).\newline
\noindent \lbrack Kuo] H. H. Kuo, Gaussian measures in Banach spaces. 
\textsl{Lecture Notes in Mathematics} \textbf{463}, Springer, 1975.\newline
\noindent \lbrack LedT] M. Ledoux and M. Talagrand, \textsl{Probability in
Banach spaces}. Ergeb. Math. (3) \textbf{33}, Springer, 1991 (paperback
2011).\newline
\noindent \lbrack LePM] R. D. LePage, V. Mandrekar, Equivalence-singularity
dichotomies from zero-one laws. \textsl{Proc. Amer. Math. Soc.} \textbf{31}
(1972), 251--254.\newline
\noindent \lbrack Lif] M. A. Lifshits, \textsl{Gaussian random functions}.
Mathematics and its Applications \textbf{322}, Kluwer, 1995.\newline
{\normalsize \noindent }[LiuR] T. S. Liu, A. van Rooij, Transformation
groups and absolutely continuous measures, \textsl{lndag. Math.} \textbf{71}
(1968), 225-231.\newline
{\normalsize \noindent }[LiuRW] T. S. Liu, A. van Rooij, J-K Wang,
Transformation groups and absolutely continuous measures II, \textsl{Indag.
Math.} \textbf{73 }(1970), 57--61.\newline
{\normalsize \noindent }[Loe] K. L\"{o}wner (=Charles Loewner), Grundz\"{u}%
ge einer Inhaltslehre im Hilbertschen Raume. \textsl{Ann. of Math.} \textbf{%
40} (1939), 816--833 (reprinted in \textsl{Charles Loewner, Collected Papers}
(ed. Lipman Bers), Birkh\"{a}user, 1988, 106-123).\newline
{\normalsize \noindent }[Lud] S. V. Ludkovsky, Properties of quasi-invariant
measures on topological groups and associated algebras. \textsl{Ann. Math.
Blaise Pascal} \textbf{6 }(1999), 33--45.\newline
{\normalsize \noindent }[Mac] G. W. Mackey, Borel structure in groups and
their duals. \textsl{Trans. Amer. Math. Soc.} \textbf{85} (1957), 134--165.%
\newline
\noindent \lbrack MarR] M. B. Marcus, J. Rosen, \textsl{Markov processes,
Gaussian processes, and local times.} Cambridge Studies in Advanced
Mathematics \textbf{100}, 2006.\newline
\noindent \lbrack McC1] M. McCrudden, On the supports of absolutely
continuous Gauss measures on connected Lie groups. \textsl{Monatsh. Math.} 
\textbf{98} (1984), 295--310.\newline
\noindent \lbrack McC2] M. McCrudden, On the supports of Gauss measures on
algebraic groups. \textsl{Math. Proc. Cambridge Philos. Soc.} \textbf{96}
(1984), 437--445.\newline
\noindent \lbrack McCW] M. McCrudden, R.M. Wood, On the support of
absolutely continuous Gauss measures.\textsl{\ Probability measures on
groups, VII (Oberwolfach, 1983)}, 379--397, \textsl{Lecture Notes in Math.} 
\textbf{1064}, Springer, 1984.\newline
\noindent \lbrack MilO] H. I. Miller and A. J. Ostaszewski, Group actions
and shift-compactness. \textsl{J. Math. Anal. Appl.} \textbf{392} (2012),
23-39.\newline
\noindent \lbrack Mor] S. A. Morris,\textsl{\ Pontryagin duality and the
structure of locally compact abelian groups. }London Math. Soc. Lecture Note
Series \textbf{29}. Cambridge University Press, 1977.\newline
\noindent \lbrack Mos] Y. V. Mospan, A converse to a theorem of Steinhaus. 
\textsl{Real An. Exch.} \textbf{31} (2005), 291-294.{\normalsize \newline
}\noindent \lbrack Neu1] J. von Neumann, Die Einf\"{u}hrung analytischer
Parameter in topologischen Gruppen. \textsl{Ann. Math.} \textbf{34} (1933),
170-179 (\textsl{Collected Works II, }366-386, Pergamon, 1961).{\normalsize 
\newline
}\noindent \lbrack Neu2] J. von Neumann, Review of [Loe], \textsl{Math.
Reviews} 1-48; MathSciNet MR0000285.{\normalsize \newline
}\noindent \lbrack Neu3] J. von Neumann, Zum Haarsche Mass in topologischen
Gruppen, \textsl{Comp. Math.} \textbf{1} (1934), 106-114 (\textsl{Collected
Works II, }445-453, Pergamon, 1961).{\normalsize \newline
}\noindent \lbrack Neu4] J. von Neumann, The uniqueness of Haar's measure. 
\textsl{Mat. Sbornik} \textbf{1} (1936), 721-734 (\textsl{Collected Works
IV, }91-104, Pergamon, 1962).{\normalsize \newline
}\noindent \lbrack Ost1] A. J. Ostaszewski, Beyond Lebesgue and Baire III:
Steinhaus's Theorem and its descendants, \textsl{Topology and its
Applications} \textbf{160} (2013), 1144-1154.\newline
\noindent \lbrack Ost2] A. J. Ostaszewski, Effros, Baire, Steinhaus and
non-separability. \textsl{Topology and its Applications} \textbf{195} (M.E.
Rudin memorial issue) (2015), 265-274.\newline
\noindent \lbrack Oxt1] J. C. Oxtoby,{\normalsize \ }Invariant measures in
groups which are not locally compact.{\normalsize \ \textsl{Trans. Amer.
Math. Soc.} }\textbf{60}{\normalsize \ }(1946), 215--237.\newline
\noindent \lbrack Oxt2] J. C. Oxtoby,{\normalsize \ \textsl{Measure and
category}, }Graduate Texts in Math. \textbf{2}, Springer, 2$^{\text{nd}}$
ed. 1980 (1$^{\text{st}}$ ed. 1972).{\normalsize \newline
}\noindent \lbrack Oxt3] J. C. Oxtoby, A commentary on [49], [62] and [63].
Pages 379-383 in [Kak3, Vol. 2].\newline
\noindent \lbrack Pan] G. R. Pantsulaia, On an invariant Borel measure in
Hilbert space. \textsl{Bull. Pol. Acad. Sci. Math.} \textbf{52} (2004),
47--51.\newline
\noindent \lbrack Par] K. R. Parthasarathy, \textsl{Probability measures on
metric spaces,} Acad. Press, 1967 (reprinted, AMS, 2005).\newline
\noindent \lbrack Pat] A. L. T. Paterson,{\normalsize \ \textsl{Amenability.}
}Math. Surveys and Mon. \textbf{29}, Amer. Math. Soc., 1988.{\normalsize 
\newline
}\noindent \lbrack Pet] {B. J. Pettis, {On continuity and openness of
homomorphisms in topological groups.} \textsl{Ann. of Math. }(2) \textbf{52}
(1950), 293--308.}\newline
\noindent \lbrack Pic] {\ S. Piccard, {Sur les ensembles de distances des
ensembles de points d'un espace Euclidien.\ }\textsl{M\'{e}m. Univ. Neuch%
\^{a}tel} \textbf{13}, 212 pp. 1939.}\newline
\noindent \lbrack Pro] V. Prokaj, A characterization of singular measures, 
\textsl{Real Anal. Exchange} \textbf{29} (2003/2004), 805--812.\newline
\noindent \lbrack Pug] O. V. Pugach\"{e}v, Quasi-invariance of Poisson
distributions with respect to transformations of configurations. \textsl{%
Dokl. Math.} \textbf{77} (2008), 420--423.\newline
\noindent \lbrack Rog] C. A. Rogers, J. Jayne, C. Dellacherie, F. Tops\o e,
J. Hoffmann-J\o rgensen, D. A. Martin, A. S. Kechris, A. H. Stone, \textsl{%
Analytic sets,} Academic Press, 1980.\newline
\noindent \lbrack Rud1] W. Rudin, \textsl{Fourier analysis on groups}.
Wiley, 1962 (reprinted 1990).\newline
\noindent \lbrack Rud2] W. Rudin, \textsl{Functional analysis}, 2$^{\text{ed}%
}$ ed., McGraw-Hill, 1991 (1$^{\text{st}}$ ed. 1973).\newline
\noindent \lbrack Sad1] G. Sadasue, On absolute continuity of the Gibbs
measure under translations.\textsl{\ J. Math. Kyoto Univ.} \textbf{41}
(2001), 257--276.\newline
\noindent \lbrack Sad2] G. Sadasue, On quasi-invariance of infinite product
measures. \textsl{J. Theoret. Probab.} \textbf{21} (2008), no. 3, 571--585.%
\newline
\noindent \lbrack Sch] L. Schwartz, \textsl{Radon measures on arbitrary
topological spaces and cylindrical measures.} Tata Institute of Fundamental
Research Studies in Mathematics \textbf{6}, Oxford University Press, 1973.%
\newline
\noindent \lbrack She] L. A. Shepp, Distinguishing a sequence of random
variables from a translate of itself. \textsl{Ann. Math. Statist.} \textbf{36%
} (1965), 1107--1112.\newline
\noindent \lbrack Shi] H. Shimomura, Quasi-invariant measures on the group
of diffeomorphisms and smooth vectors of unitary representations. \textsl{J.
Funct. Anal.} \textbf{187} (2001), 406--441.\newline
\noindent \lbrack Sim] S. M. Simmons, A converse Steinhaus theorem for
locally compact groups. \textsl{Proc. Amer. Math. Soc.} \textbf{49} (1975),
383-386.{\normalsize \newline
}\noindent \lbrack Sko] A. V. Skorohod, \textsl{Integration in Hilbert space.%
} Ergebnisse Math. \textbf{79}, Springer, 1974.\newline
\noindent \lbrack Smo] W. Smole\'{n}ski, On quasi-invariance of product
measures. \textsl{Demonstratio Math.} \textbf{11} (1978), no. 3, 801--805.%
\newline
\noindent \lbrack Sol] S. Solecki, Amenability, free subgroups, and Haar
null sets in non-locally compact groups. \textsl{Proc. London Math. Soc.} 
\textbf{(3) 93} (2006), 693--722.\newline
\noindent \lbrack Ste] H. Steinhaus, Sur les distances des points de mesure
positive.\textsl{\ Fund. Math.} \textbf{1} (1920), 83-104.{\normalsize 
\newline
}\noindent \lbrack Str] D. W. Stroock, \textsl{Probability theory. An
analytic view.} 2$^{\text{nd}}$ ed. Cambridge University Press, 2011 (1$^{%
\text{st}}$ ed. 1993).{\normalsize \newline
\noindent }[Tar] V. Tarieladze, Characteristic functionals of probabilistic
measures in DS-groups and related topics. \textsl{J. Math. Sci.} \textbf{211}
(2015), 137--296.{\normalsize \newline
\noindent }[TopH] F. Tops\o e and J. Hoffmann-J\o rgensen, Analytic spaces
and their application, in: [Rog].{\normalsize \newline
\noindent }[Wei] A. Weil,{\normalsize \textsl{\ L'int\'{e}gration dans les
groupes topologiques}, }Actualit\'{e}s Scientifiques et Industrielles 1145,
Hermann, 1965 (1$^{\text{st }}$ ed. 1940).\newline
{\normalsize \noindent }[Xia] D. X. Xia,\textsl{\ Measure and integration
theory on infinite-dimensional spaces. Abstract harmonic analysis.} Pure and
App. Math. \textbf{48}. Academic Press, 1972.\newline
\noindent \lbrack Yam1] Y. Yamasaki, Translationally invariant measure on
the infinite-dimensional vector space. \textsl{Publ. Res. Inst. Math. Sci. }%
\textbf{16} (1980), 693--720.\newline
\noindent \lbrack Yam2] Y. Yamasaki, \textsl{Measures on
infinite-dimensional spaces.} World Scientific, 1985.\newline
\noindent \lbrack Yos] K. Yosida, \textsl{Functional analysis}. Reprint of
the sixth (1980) edition. Classics in Mathematics. Springer, 1995.

\bigskip

\noindent Mathematics Department, Imperial College, London SW7 2AZ;
n.bingham@ic.ac.uk \newline
Mathematics Department, London School of Economics, Houghton Street, London
WC2A 2AE; A.J.Ostaszewski@lse.ac.uk

\end{document}